\documentclass[11pt,a4paper]{article}
\usepackage{epsfig}

\usepackage{amsmath}
\usepackage{amssymb}
\usepackage{amsthm}
\usepackage{amsfonts}
\usepackage{graphicx}
\usepackage{epic}
\usepackage{eepic}

\usepackage{hyperref}

\usepackage{color}
\usepackage{bbm,verbatim}

\def\E{{\mathbb E}}
\def\R{{\mathbb R}}
\def\P{{\mathbb P}}
\def\Z{{\mathbb Z}}
\def\N{{\mathbb N}}

\def\T{{\mathbb T}}

\def\calP{{\mathcal P}}

\def\Span{\rm span}
\def\eps{\epsilon}
\def\Tg{\mathbb{T}}

\def\Spec{\mathcal{S}}

\def\md{\mid}
\def\Bb#1#2{{\def\md{\bigm| }#1\bigl[#2\bigr]}}

\def\Pb{\Bb\P}
\def\Eb{\Bb\E}

\newtheorem{theorem}{Theorem}[section]
\newtheorem{lemma}[theorem]{Lemma}
\newtheorem{corollary}[theorem]{Corollary}
\newtheorem{proposition}[theorem]{Proposition}
\newtheorem{definition}[theorem]{Definition}
\newtheorem{example}[theorem]{Example}

\newtheorem{conjecture}[theorem]{Conjecture}
\newtheorem{question}[theorem]{Question}

\author{Erik I. Broman \and Christophe Garban \and Jeffrey E. Steif}

\date{\today}

\title{Exclusion Sensitivity of Boolean Functions}

\begin{document}
  \maketitle

  \begin{abstract}
Recently the study of noise sensitivity and noise stability of Boolean
functions has received considerable attention. The purpose of this paper is to
extend these notions in
a natural way to a different class of perturbations, namely those arising
from running the symmetric exclusion process for a short amount of time.
In this study, the case of monotone
Boolean functions will turn out to be of particular interest. We show that
for this class of functions, ordinary noise sensitivity and noise
sensitivity with respect to the complete graph exclusion process are equivalent.
We also show this equivalence with respect to stability.

After obtaining these fairly general results, we study
``exclusion sensitivity'' of critical percolation in more detail with respect to
medium-range dynamics.
  The exclusion dynamics, due to its conservative nature, is in some sense more physical than the
classical i.i.d. dynamics. Interestingly, we will see
  that in order to obtain a precise understanding of the exclusion sensitivity of percolation,
we will need to describe how typical spectral sets of percolation
  diffuse under the underlying exclusion process.




  \end{abstract}

\newpage

\tableofcontents

\begin{figure}\label{f.Z2site}
\begin{center}
\includegraphics[width=\textwidth]{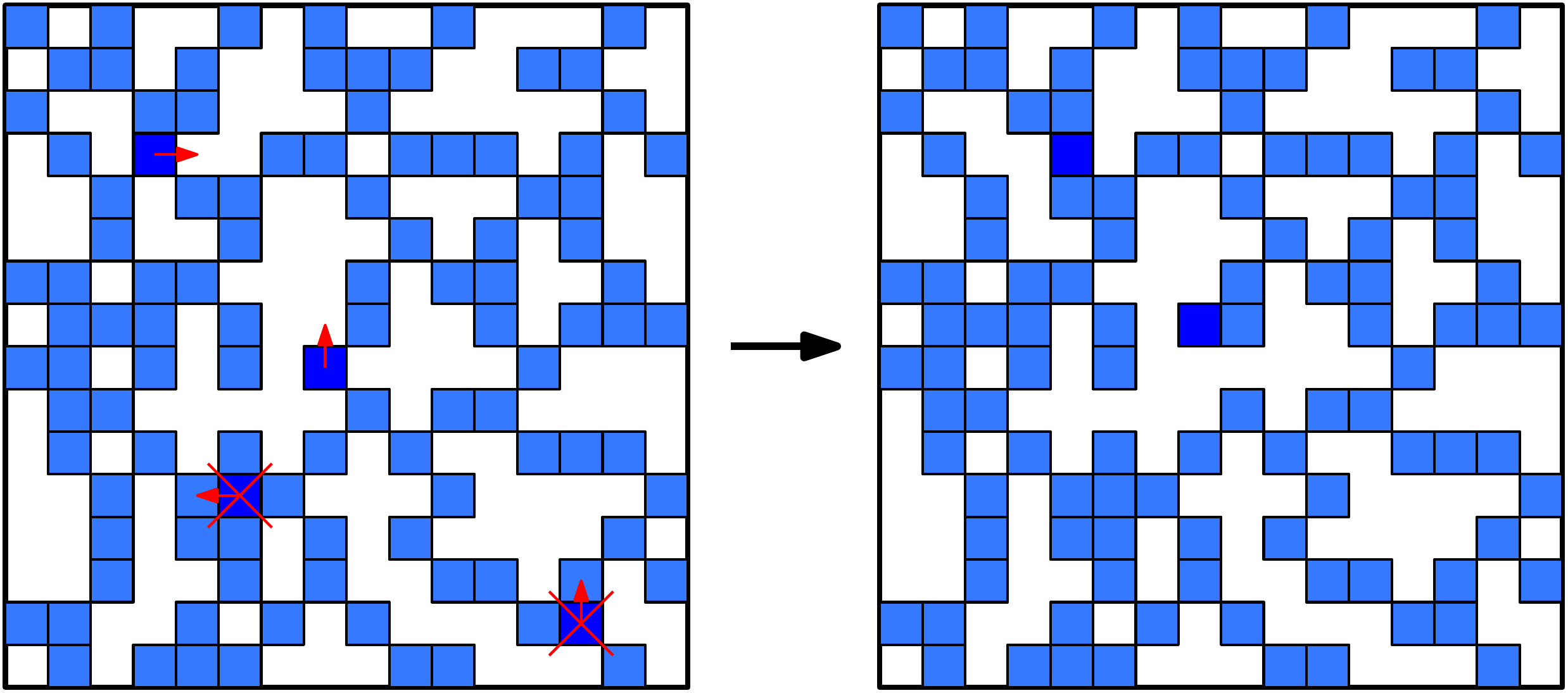}
\end{center}
\caption{\bf This picture illustrates the nearest-neighbor dynamics for {\bf site}-percolation on 
$\Z^2$.}\label{f.dynamics}
\end{figure}

\section{Introduction}\label{secIntro}
The concept of noise sensitivity was introduced in \cite{BKS} and is defined as follows.
Let $V_n$ be some finite set and $\omega\in\{0,1\}^{V_n}$ be such that
$\omega_i=1$ with probability $1/2$
independently for every $i\in V_n$. Given $\epsilon \in (0,1)$, resample each
$\omega_i$ independently with probability $\epsilon$ and call the resulting configuration
$\omega^{\epsilon}$.
Throughout, $\{V_n\}_{n\in\N}$ will be an increasing sequence of finite sets.
Given such a sequence we will consider
sequences of functions $\{f_{n}\}_{n \in \N}$ such that
$f_{n}:\{0,1\}^{V_n}\rightarrow \{-1,1\}$. These functions will
sometimes be referred to as \emph{Boolean functions}. We have the
following central definition introduced in \cite{BKS}.
\begin{definition} \label{defnoise1}
The sequence of functions $f_{n}:\{0,1\}^{V_n}\rightarrow \{-1,1\}$
is noise sensitive (NS) if for any $\epsilon>0$,
\begin{equation} \label{eqn:NS}
\lim_{n\rightarrow \infty} \E[f_{n}(\omega)f_{n}(\omega^{\epsilon})]-\E[f_{n}(\omega)]^2=0.
\end{equation}
\end{definition}

\noindent
{\bf Remark:}  It is known that if the limit above is 0 for some
$\epsilon\in (0,1)$,
then it is 0 for all such $\epsilon$; this will be clear from Theorem
\ref{thmBKS} below. We also mention that it is the case that
$\E[f_{n}(\omega)f_{n}(\omega^{\epsilon})]-\E[f_{n}(\omega)]^2$ is
nonnegative and decreasing in $\epsilon$; this is seen via a
Fourier representation, see Section \ref{secFourier}.

\bigskip

Closely related to NS is the notion of noise \emph{stability} defined as
follows.

\begin{definition} \label{defstab1}
The sequence of functions $f_{n}:\{0,1\}^{V_n}\rightarrow \{-1,1\}$
is noise stable (NStable) if for any $\delta>0$ there exists an $\epsilon>0$ such that
for every $n$,
\[
\P[f_{n}(\omega)\neq f_{n}(\omega^{\epsilon})]\le \delta.
\]
\end{definition}

\noindent
{\bf Remark:} In some sense, noise stability is the opposite of noise sensitivity. It is however not
very difficult to find examples that are neither noise sensitive nor noise stable.
\medskip

The concepts that we will study in this paper are intimately related
to noise sensitivity and noise stability. We first describe noise sensitivity and noise stability
in slightly different terms.
Let $\pi_{p}$ denote product measure on $\{0,1\}^{V_n}$ with density $p$,
and consider the following process on $\{0,1\}^{V_n}$.
Pick $\omega_0$ according to $\pi_{1/2}$ and associate to every $k\in
V_n$ a Poisson clock
with rate 1. Whenever the Poisson clock rings for some $k\in V_n$, the
site $k$ changes its value to $0$ or $1$ with equal probability independently
of all past decisions. This results in a stochastic
process $\{\omega_t\}_{t\geq 0}$.
Note that for fixed $t$ the probability that $k\in V_n$ has been
updated by time $t$ is $1-e^{-t}$
independently of every other vertex. Therefore,
the pair $(\omega,\omega^{\epsilon})$ introduced earlier
will have the same distribution as
$(\omega_0,\omega_t)$ where $\epsilon=1-e^{-t}$.
It follows that
the sequence of functions $f_{n}:\{0,1\}^{V_n}\rightarrow \{-1,1\}$
is NS if and only if for some (equivalently for all) $t>0$,
\begin{equation} \label{eqnNSequiv}
\lim_{n\rightarrow \infty} \E[f_{n}(\omega_0)f_{n}(\omega_{t})]-\E[f_{n}(\omega_0)]^2=0.
\end{equation}
Furthermore, $\{f_n\}_{n \in \N}$ is NStable if and only if for every $\delta>0$ there exists a $t>0$ such that
for every $n$,
\begin{equation} \label{eqnNStableequiv}
\P[f_{n}(\omega_0)\neq f_{n}(\omega_t)]\le \delta.
\end{equation}
The process $\{\omega_t\}_{t\geq 0}$ will be referred to as the
Simple Non-interacting Particle System or SNPS for short.

Equations (\ref{eqnNSequiv}) and (\ref{eqnNStableequiv}) are the starting
points of what we do below.
Let $\{G_n\}=\{(V_n,E_n)\}$ be an increasing sequence of finite graphs and
$\{\alpha_n\}_{n\in\N}$ be a sequence of positive real numbers.
In this paper we will consider
the {\em symmetric exclusion process} on $\{0,1\}^{V_n}$
with rate $\alpha_n$.
It is natural to work with the symmetric
exclusion process since the (unique) invariant measure for SNPS is also
an invariant measure for the symmetric exclusion process (although there
are many other invariant measures for this latter process).
The symmetric exclusion process can be defined
in the following way (see \cite{SIS} for a survey on the exclusion process).

The first step is to define a Markov process $\pi^n_t$ on the set of
permutations of $V_n$ as follows.
For every edge $e\in E_n$, associate an independent Poisson clock with parameter $\alpha_n$.
When this Poisson clock rings, we interchange the endpoints  of $e$.
The permutation at time $t$, $\pi^n_t$, is obtained by composing all of the
transpositions that have occurred up to time $t$
in the order in which they occurred.

For $S,S'\subseteq V_n$, we let
$$
\P^n_t(S,S')=\P(\pi^n_t(S)=S').
$$
Note that this probability is 0 if $|S|\neq |S'|$.
It is not hard to see that $\pi^n_t$ and $(\pi^n_t)^{-1}$ have the same
distribution and hence that $\P^n_t$ is a symmetric matrix. Furthermore we let
$S_t:=\pi^n_t(S)$.

The symmetric exclusion process on $(G_n,\alpha_n)$ starting from $\eta^n_0$ is
the process $\{\eta^n_t\}$ taking values in $\{0,1\}^{V_n}$ defined by
$$
\eta^n_t(x)=\eta^n_0((\pi^n_t)^{-1}(x)).
$$

We will always take $\eta^n_0\in\{0,1\}^{V_n}$ to be chosen according to
$\pi_{1/2}$. It is then clear that $\{\eta^n_t\}$ is a Markov process
defined in $\{0,1\}^{V_n}$ which evolves
as follows. If $t$ is a time when the clock associated to
the edge $e=(v_1,v_2)$ rings we let $\eta^n_{t}(v_1)=\eta^n_{t^-}(v_2)$ and
$\eta^n_{t}(v_2)=\eta^n_{t^-}(v_1)$.
For convenience we will from now on refer to the symmetric exclusion process
simply as the ``exclusion process''.

As a more philosophical remark,
the mathematical model that we are studying would be relevant if we want to model
some type of conservative dynamics. For example, we might want to understand
what happens if we 'shake our system' a little, moving the 1's and 0's
around.

\medskip\noindent
{\bf Fact:} \\
1. It is elementary to check, due to a detailed balance condition,
that the set of homogeneous product measures are stationary distributions
for the exclusion process. In addition, since the
Markov chain is not irreducible (the irreducible classes are determined by
specifying the number of ones in each component of the graph), there are other
stationary distributions as well. Here and later ``particle conservation''
simply refers to the fact that the number of sites $v$ with $\eta^n_t(v)=1$ is
constant in time.\\
2. One can also discuss nonsymmetric exclusion processes but these do not
have a representation in terms of a random permutation as above.

\medskip\noindent
Given $G_n=(V_n,E_n)$ and $\alpha_n$ as above,
which then yields the process $\{\eta^n_t\}_{t\geq 0}$,
and given $f_n:\{0,1\}^{V_n} \rightarrow \{-1,1\}$, let
\[
N(f_n,t,G_n,\alpha_n):={\mathbb E}[f_n(\eta^n_0)f_n(\eta^n_{t})]-{\mathbb E}[f_n(\eta^n_0)]^2.
\]
The following fact will be useful; it will be established in Section
\ref{secXSusingFourier}.

\begin{proposition} \label{propalltsame}
Given $\{G_n,\alpha_n\}_{n \in \N}$ and $\{f_n\}_{n \in \N}$ as above,
then, for any $t,s >0$,
$\lim_{n\rightarrow \infty} N(f_n,t,G_n,\alpha_n)=0$ if and only if
$\lim_{n\rightarrow \infty} N(f_n,s,G_n,\alpha_n)=0$.
\end{proposition}

\noindent
The following definition is analogous to Definition \ref{defnoise1}.

\begin{definition}\label{defXS}
The sequence of functions
$f_n:\{0,1\}^{V_n} \rightarrow \{-1,1\}$
is XS (eXclusion Sensitive) with respect to $\{G_n,\alpha_n\}_{n\in\N}$
if for some $t>0$ (equivalently, by Proposition \ref{propalltsame}, for all $t>0$)
\[
\lim_{n\rightarrow \infty} N(f_n,t,G_n,\alpha_n)=0.
\]
\end{definition}

\medskip\noindent
{\bf Remark:} A natural generalization of the above concept is
to allow $\alpha_n(\{v_1,v_2\})$ to depend on the 2 element
subset $\{v_1,v_2\}\subseteq V_n$. Obviously in the above, we are just
considering the case where $\alpha_n(\{v_1,v_2\})=0$
for every $\{v_1,v_2\} \not \in E_n$ and
$\alpha_n(\{v_1,v_2\})=\alpha_n$ for every $\{v_1,v_2\} \in E_n$.
It is also possible to generalize
the concept of sensitivity to a range of other interacting particle systems.
However, since we only deal with the exclusion process
here, we don't give these obvious definitions.
\medskip

Analogously to noise stability, we have the following definition for the
exclusion process.

\begin{definition}\label{defXStable}
The sequence of functions $f_n:\{0,1\}^{V_n} \rightarrow \{-1,1\}$
is XStable (eXclusion Stable)
with respect to $\{G_n,\alpha_n\}_{n\in\N}$ if for any $\delta>0$,
there exists a $t>0$ such that for every $n$,
\[
\P[f_n(\eta^n_0)\neq f_n(\eta^n_{t})]\le \delta.
\]
\end{definition}

\noindent
{\bf Remark:} A similar remark to the one following Definition \ref{defstab1}
applies also in this case.

\medskip

We want to warn the reader that in some cases, $f_{n}$ will be
defined in terms of a specific graph structure which may be different
from the graph $G_n$. For example we will study functions $f_{n}$ defined
naturally in terms of the nearest neighbor graph on $\{1,\ldots,n\}^2$,
while at the same time considering the complete graph dynamics (defined below).

Throughout the paper we will consider various choices of
$V_n,E_n$ and $\alpha_n$. We make the following assumption throughout
the rest of the paper.

\medskip\noindent
{\bf Assumption:} For each $n$, $\alpha_n\le 1/d(G_n)$ where
$d(G_n)$ is the maximum degree in $G_n$.

\medskip\noindent
A particularly important example will be when $G_n$ is the complete
graph on $n$ vertices and where we let $\alpha_n=1/|V_n|$.

\begin{definition}\label{defCGXS}
The sequence of functions $f_n:\{0,1\}^{V_n} \rightarrow \{-1,1\}$ is
CGXS (Complete Graph eXclusion Sensitive)
if they are XS with respect to the complete graph and $\alpha_n=1/|V_n|$.
\end{definition}

\noindent
Analogously we have the following definition.

\begin{definition}\label{defCGXStable}
The sequence of functions $f_n:\{0,1\}^{V_n} \rightarrow \{-1,1\}$ is
CGXStable if they are XStable with respect to the complete graph and $\alpha_n=1/|V_n|$.
\end{definition}

We will often let $|V_n|=n$, but sometimes other choices are more natural;
for example, when $G_{n}=\Lambda^d_{n}$ where $\Lambda^d_{n} \subseteq
\Z^d$ is the box of side length $n$, it is more natural to let
$V_n=\Lambda^d_{n}$ and $E_n$ be the set of edges induced by $\Z^d$.
In this case we choose $\alpha_n=1/(2d)$.

\begin{definition}\label{defNNXS}
The sequence of functions
$f_{n}:\{0,1\}^{\Lambda^d_{n}} \rightarrow \{-1,1\}$ is d-NNXS (Nearest
Neighbor eXclusion Sensitive in dimension $d$) if it is XS with respect
to the sequence of graphs $\Lambda_{n}^d$ and $\alpha_n=1/(2d)$.
\end{definition}

\noindent
Analogously we have the following definition.

\begin{definition}\label{defNNXStable}
The sequence of functions $f_n:\{0,1\}^{\Lambda^d_{n}} \rightarrow \{-1,1\}$ is
d-NNXStable if it is XStable with respect to the sequence of graphs
$\Lambda_{n}^d$ and $\alpha_n=1/(2d)$.
\end{definition}

The following example distinguishes between 1-d nearest neighbor dynamics
and the complete graph dynamics.

\begin{example} \label{ex7}
Let $f_n$ be $\sum_{i=1}^{n/2} \omega_i ({\rm mod } \, 2)$ where
$V_n=[n]:=\{1,2,\ldots,n\}$ (also known as the parity
function on the first $n/2$ bits). It is easy to show that it is CGXS
but not 1-NNXS; in fact it is 1-NNXStable. We leave the details to the reader.
\end{example}

\medskip\noindent
Example \ref{ex7} is somewhat unsatisfactory since in $d=1$ it is easy to
keep different parts of $[n]$ ``blocked off'' from each other. Later on,
in Section \ref{prop:CGXSNOTIMPLY2-NNXSSECTION},
we give some examples which are CGXS but not 2-NNXS.

These examples partly motivate the following questions.
\begin{question} \label{Q1}
For a sequence $\{G_n,\alpha_n\}_{n\in\N}$
and a sequence of functions $f_{n}:\{0,1\}^{V_n}\rightarrow \{-1,1\}$,
is it the case that
\[
\{f_{n}\}_{n\in \N} \textrm{ is XS } \Rightarrow \{f_{n}\}_{n\in \N} \textrm{ is CGXS }?
\]
Is it the case that
\[
\{f_{n}\}_{n\in \N} \textrm{ is CGXStable } \Rightarrow \{f_{n}\}_{n\in \N} \textrm{ is XStable } ?
\]
\end{question}

\medskip\noindent
{\bf Remark:}
At this point, we do not have an answer to these questions.
More generally, one could ask whether one particular sequence of graphs
is ``more sensitive'' (or ``more stable'')
than another particular sequence of graphs; of course, this question
is quite vague.
\medskip

One might ask how the notion of NS and NStable relate to various variants of XS
and XStable (with respect to some sequence $\{G_n,\alpha_n\}_{n\in\N}$).
Our next result gives a partial answer.

\begin{theorem} \label{thm2}
1. If, for some choice of sequence $\{G_n,\alpha_n\}_{n\in\N}$,
the sequence of functions $f_{n}:\{0,1\}^{V_n}\rightarrow \{-1,1\}$
is XS with respect to $\{G_n,\alpha_n\}_{n\in\N}$,
then the sequence is also NS.\\
2. If the sequence $\{f_n\}_{n \in \N}$ is NStable, then the sequence is XStable
with respect to any $\{G_n,\alpha_n\}_{n\in\N}$.
\end{theorem}

\medskip\noindent
{\bf Remark:}
Note that the notion of NS has no $\alpha_n$-dependence and that this holds due to our earlier
"Assumption".

Note further that both converses are trivially false.
Letting $f_n$ be $\sum_i \omega_i ({\rm mod}\, 2)$ (i.e., the parity
function), it is easy to see that $\{f_n\}_{n \in \N}$ is NS while
it follows trivially from particle conservation that $\{f_n\}_{n \in \N}$
is not XS for any sequence $\{G_n,\alpha_n\}_{n\in\N}$; in fact, $\{f_n\}_{n \in \N}$ is
clearly XStable for every sequence $\{G_n,\alpha_n\}_{n\in\N}$.

\medskip\noindent
Given  a sequence of events ${\cal A}_n\subseteq \{0,1\}^{V_n}$, we let
$f_n=2I_{{\cal A}_n}-1$,
where $I_{{\cal A}_n}$ is the indicator function of the event ${\cal  A}_n$.
We say that the sequence of events $\{{\cal A}_n\}_{n\in \N}$
is XS (XStable) with respect to $\{G_n,\alpha_n\}_{n\in \N}$
if the sequence $\{f_n\}_{n\in \N}$ is XS (XStable)
with respect to $\{G_n,\alpha_n\}_{n\in \N}$.

\medskip\noindent
A central concept which arises in this area is the notion of pivotality and
influence. Letting $\omega_i$ be $\omega$ with the $i$th bit flipped,
we say that $i$ is \emph{pivotal} for $f$ if $f(\omega)\neq f(\omega_i)$
and we define the \emph{influence} of the variable $i$ for a Boolean function
$f$ to be
$$
I_i(f):=\P(f(\omega)\neq f(\omega_i)).
$$
Letting
\[
I\!I(f):=\sum_{i}I_i(f)^2,
\]
one of the main results in \cite{BKS} is the following.

\begin{theorem} \label{BKSMAIN} \cite{BKS}
If a sequence of Boolean functions $\{f_{n}\}_{n\in \N}$ satisfies
$I\!I(f_n)\to 0$  as $n\to\infty$, then the sequence is noise
sensitive.
\end{theorem}

One of our main results is that this condition is also sufficient for
complete graph exclusion sensitivity.

\begin{theorem} \label{thmmonNSimpliesCGXS}
If $\{f_n\}_{n\in \N}$ is a sequence that satisfies
$\lim_{n \rightarrow \infty}I\!I(f_n)=0$, then the sequence is CGXS.
\end{theorem}

The parity function shows that the converse of Theorem \ref{BKSMAIN} is false.
However the converse is true (see \cite{BKS}) for the
very natural and often occurring class of \emph{monotone} Boolean functions.
A Boolean function $f_n$
is \emph{monotone} if $f_n(\omega)\leq f_n(\omega')$ for every $\omega,\omega'$
such that $\omega_i\leq \omega'_i$ for every $i$.
This together with Theorem \ref{thmmonNSimpliesCGXS} yields the following corollary.

\begin{corollary} \label{cor:monotone}
If $\{f_n\}_{n\in \N}$ are monotone and NS, then the sequence is CGXS.
\end{corollary}

\medskip\noindent
{\bf Remark:}
In combination with Theorem \ref{thm2}, this shows that for monotone
functions NS is equivalent to CGXS.

\medskip\noindent
It turns out that not only are CGXS and NS equivalent for monotone functions,
but CGXStable and NStable are also equivalent for monotone functions.
One direction is given by the second part of Theorem \ref{thm2}
while the other direction is given next.

\begin{theorem} \label{thmmonCGXStableimpliesNstable}
If $\{f_n\}_{n\in \N}$ are monotone and
CGXStable, then the sequence is NStable.
\end{theorem}

\medskip\noindent
The method of proof for Theorem \ref{thmmonNSimpliesCGXS} relates
\emph{quantitative} noise sensitivity to \emph{quantitative} complete graph exclusion sensitivity.
We must first give a little background. In the definition of NS, $\epsilon$ is fixed.
However, in the study of quantitative noise sensitivity, one allows
$\epsilon$ to depend on $n$ and therefore use the notation $\epsilon_n$.
One then asks how quickly
$\epsilon_n$ can go to 0 but still yielding (\ref{eqn:NS}). Considering $t_n$ in the obvious
analogous way we have the following quantitative variant of Theorem \ref{thmmonNSimpliesCGXS}
which follows from its proof.
\begin{corollary} \label{corr1}
Consider the sequences $\{\eps_n\}_{n \in \N}$ and $\{t_n\}_{n \in \N}$ where $\eps_n=1-e^{-t_n}$ for every $n$
and let the sequence $\{f_n\}_{n \in \N}$ satisfy $\lim_{n \to \infty}I\!I(f_n)=0$.
Then the sequence $\{f_n\}_{n \in \N}$
is NS with respect to $\{\eps_n\}_{n \in \N}$ iff it is CGXS with respect to $\{t_n\}_{n \in \N}$.
\end{corollary}

The sequence of Boolean functions which have been most studied in the theory
of noise sensitivity are percolation crossings which are defined as follows.
We consider critical percolation on the triangular lattice $\T$ where
each vertex is 1 (open) or 0 (closed) with probability $1/2$ each.
Fix $a,b>0$, by $[0,an]\times[0,bn]$, we mean the subset of
$\T$ which is contained in the corresponding rectangle in Euclidean space.
We let $C_{n}$ denote the event that there exists a left-to-right crossing
by open sites (i.e. sites with value 1) of $[0,an]\times[0,bn]$.
Noise sensitivity of
this event has been extensively studied; see \cite{BKS}, \cite{SS} and
\cite{GPS}.  In \cite{BKS}, noise sensitivity of this event
was established as stated next; in \cite{SS} and \cite{GPS},
quantitative versions were established.

\begin{theorem} \label{thmPCNS} \cite{BKS}
The sequence $\{C_{n}\}_{n\in \N}$ is NS.
\end{theorem}

Using information concerning the geometrical structure of the
\emph{spectral sample}
(see Section \ref{secFourier} for definitions) which has been established
in \cite{GPS}, we will prove the following result.

\begin{theorem} \label{thm5}
The sequence $\{C_{n}\}_{n\in \N}$ is CGXS.
\end{theorem}

\medskip\noindent
Note that, due to Theorem \ref{thmPCNS}, Theorem \ref{thm5} is a special case of
Corollary \ref{cor:monotone}.  However, we
reprove Theorem \ref{thm5} since it demonstrates a different and
interesting method of proof based on spectral techniques
which also can be used to handle a more ``realistic'' medium-range
dynamics in the case of critical percolation on the triangular graph.
In addition, the proof does not use the monotonicity property, and is therefore potentially
applicable to non-monotone cases.

Ideally, one would like to consider a nearest-neighbor dynamics, but as we will explain
in Section \ref{secdiscussions}, a spectral approach in this case seems highly non-trivial.
See Figure \ref{f.dynamics}.

The medium-range dynamics is defined as follows. Let $0<\alpha<1$.
For each $n$, we define an exclusion dynamics on $\{0,1\}^{V(\Tg)}$
by transposing any $\{x,y\}$ with $\|x-y\|\le n^\alpha$ at rate $1/n^{2\alpha}$.
We denote this medium-range exclusion process by $\{\eta^{\alpha,n}_t\}$.
We have the following result.
\begin{theorem}\label{th.SR}
For any $\alpha>0$ and $t>0$ and letting $f_n=I_{C_n}$, we have
\[
\lim_{n \to \infty}\E[f_n(\eta^{\alpha,n}_0)f_n(\eta^{\alpha,n}_t)]-\E[f_n(\eta^{\alpha,n}_0)]^2 = 0.
\]
\end{theorem}

\noindent
{\bf Remarks:} 1. One could prove the same statement on the square lattice $\Z^2$.
However, much more is known about the triangular grid case and therefore we stick
to this for simplicity. \\
2. Since the statement assumes a dynamics on a full-plane configuration,
 we are not using our usual notation of being XS with respect to some
$\{G_n,\alpha_n\}_{n\in \N}$ but clearly this is completely analogous.
One could define a similar medium-range dynamics on the finite
boxes $[0,an]\times[0,bn]$,
but there would be some slightly unpleasant
boundary effects and so we choose the above setup instead.

\medskip
We next mention two other results in \cite{BKS} which
turn out to carry through for the exclusion process on the complete graph.
The first result gives us an equivalent description of noise
sensitivity which might appear stronger.

\begin{proposition} \label{BKSequiv} \cite{BKS}
A sequence of Boolean functions $\{f_{n}\}_{n\in \N}$
is noise sensitive if and only if for all $\epsilon >0$,
$$
\E[f_n(\omega^\epsilon) | \, \omega]\to \E[f_n(\omega)]
$$
in probability as $n\to\infty$.
\end{proposition}

The above is crucially used in the following interesting result from
\cite{BKS}.

\begin{theorem} \label{perccrossings} \cite{BKS}
Consider the SNPS run for time 1
and the percolation crossing events $C_n$.
If we let $S_n$ be the number of times during $[0,1]$ that
we switch from $C_n$ occurring to not occurring, then
$S_n \to\infty$ in probability as $n\to\infty$.
\end{theorem}

It turns out that one can check that the analogues of these results
hold in our new setting. Since these results can be proved in a way
very similar to the corresponding results in \cite{BKS},
we skip the proofs and simply state the results.

\begin{proposition} 
Given $\{G_n,\alpha_n\}_{n \in \N}$ and $\{f_n\}_{n \in \N}$,
the sequence is XS if and only if for all $t>0$,
$$
\E[f_n(\eta^n_t)|\, \eta^n_0]\to \E[f_n(\eta^n_0)]
$$
in probability as $n\to\infty$.
\end{proposition}

\begin{theorem} 
Consider the exclusion process on the complete graph run for time 1
and the percolation crossing events $C_n$.
If we let $S_n$ be the number of times during $[0,1]$ that
we switch from $C_n$ occurring to not occurring, then
$S_n \to\infty$ in probability as $n\to\infty$.
\end{theorem}

It turns out that the proof of Theorem \ref{thmmonNSimpliesCGXS}
also yields results for \emph{ordinary} quantitative noise sensitivity.
In \cite{MO}, it was proved that
for any $\delta >0$, there is a sequence of nondegenerate
monotone Boolean functions on $n$ variables for which
(\ref{eqn:NS}) goes to 0 when $\epsilon_n=1/n^{1/2-\delta}$.
(In fact, they actually proved the weaker fact that
$\P(f_{n}(\omega)\neq f_{n}(\omega^{\epsilon_n}))$ does not tend to 0
but their methods yield the stronger claim.) Next, it is known
(see \cite{GS} for example) that this is sharp in the sense that for
any sequence of nondegenerate monotone Boolean functions on $n$ variables, if
$\epsilon_n=c/n^{1/2}$ for any fixed $c$,
then (\ref{eqn:NS}) fails. This says that we cannot,
when sticking to nondegenerate monotone Boolean functions on $n$ variables
and using $\epsilon_n=c/n^{1/2}$ for any $c$, have $f_{n}(\omega)$ and
$f_{n}(\omega^{\epsilon_n})$ asymptotically uncorrelated.
However, interestingly, this does not necessarily imply that they are
asymptotically perfectly correlated; in \cite{MO} they construct an example, based
on earlier work of Talagrand, of a sequence of
nondegenerate monotone Boolean functions on $n$ variables,
such that when $\epsilon_n=C/n^{1/2}$ for some sufficiently large constant
$C$,
$$
\inf_n \P(f_{n}(\omega)\neq f_{n}(\omega^{\epsilon_n}))>0.
$$
The following proposition says, paradoxically, that
any such example cannot be noise sensitive (in the usual non-quantitative sense).

\begin{proposition} \label{prop:paradox}
Assume that $I\!I(f_n)\to 0$  as $n\to\infty$. Then for any constant
$C$, if we take $\epsilon_n=C/n^{1/2}$, then
$$
\lim_{n\to\infty}\P(f_{n}(\omega)\neq f_{n}(\omega^{\epsilon_n}))=0.
$$
In particular, this conclusion holds if $\{f_n\}_{n\in \N}$ is monotone and NS.
\end{proposition}

\noindent{\bf Remark:} Informally, this result shows that if one wishes to construct a sequence of
monotone Boolean functions with at least some (uniformly in $n$) spectral mass (see Section \ref{secFourier}) ``near $\sqrt{n}$'',
then one is forced to also have some mass on finite frequencies.
Again this indicates the ``rigidity'' of monotone Boolean functions.

\medskip\noindent
We mention that in \cite{BKS}, another type of noise is briefly considered.
In their model, a fixed size random subset of the bits are flipped. Clearly,
parity is no longer sensitive for such a noise but it is proved in
\cite{BKS} that, with a certain range of the size of the number of bits flipped,
if $I\!I(f_n)\to 0$, then we have NS with respect to this new type of noise.

\medskip\noindent
The rest of the paper is outlined as follows. In Section \ref{secFourier} we introduce
the Fourier-Walsh decomposition and state the standard spectral characterization of
both noise sensitivity and noise stability. In Section \ref{secXSusingFourier},
we give an analogous characterization of both XS and XStable. Proposition
\ref{propalltsame} will follow from this characterization of XS.
In Section \ref{prop:CGXSNOTIMPLY2-NNXSSECTION},
some examples which are CGXS but not 2-NNXS are presented.
The proof of Theorem \ref{thm2} is given in Section \ref{secproofs}
while the proofs of Theorem \ref{thmmonNSimpliesCGXS}, Theorem \ref{thmmonCGXStableimpliesNstable},
Corollary \ref{corr1} and Proposition \ref{prop:paradox}
are given in Section \ref{sect:monotonicity}.
The (re)proof of Theorem \ref{thm5} is given in Section \ref{secthm5}
while we prove Theorem \ref{th.SR} in Section \ref{s.shortrange}.
Examples illustrating various phenomena also appear throughout
the various sections. Finally in Section \ref{secdiscussions},
we conclude with some further discussion.

\section{Fourier-Walsh decomposition and the Spectral Sample} \label{secFourier}
For some set $V$, $\omega\in\{0,1\}^{V}$ and $i\in V$, we define
\[
\chi_i(\omega)=\left\{
\begin{array}{cc}
-1 & \textrm{if } \omega_i=1 \\
1 & \textrm{if } \omega_i=0.
\end{array}
\right.
\]
Furthermore, for $S\subseteq V,$ let $\chi_S(\omega):=\prod_{i\in
  S}\chi_i(\omega)$. (In particular, $\chi_\emptyset$ is the constant
function
1.)
The set $\{\chi_S\}_{S\subseteq V}$ forms an orthonormal basis for the set
of functions $f:\{0,1\}^{V}\mapsto \R$. We can therefore expand such functions
\[
f(\omega)=\sum_{S\subseteq V} \hat{f}(S)\chi_S(\omega),
\]
where $\hat{f}(S):=\E[f\chi_S]$. This is sometimes called the Fourier-Walsh decomposition.
Observe that when $f:\{0,1\}^{V}\mapsto \{-1,1\}$ we have that $1=\E[f^2]=\sum_{S\subseteq V} \hat{f}(S)^2$.
In this case we can define a random variable ${\cal S}$
taking values in the set of subsets of $V$ by
letting the distribution of ${\cal S}$ be given by
$\P({\cal S}=S)=\hat{f}(S)^2$ for every $S\subseteq V$.
The random variable ${\cal S}$ is called the spectral sample, and
its dependency on $f$ will sometimes be indicated by writing ${\cal S}_f$.
Furthermore, from now on, we will not stress in the notation
that $S\subseteq V$.

The following crucial formula from \cite{BKS} is easily proved.
\begin{equation} \label{keyformula}
\E[f(\omega)f(\omega^{\epsilon})]-\E[f(\omega)]^2=
\E[(1-\epsilon)^{|{\cal S}_{f}|}I_{{\cal S}_{f}\neq\emptyset}].
\end{equation}
(The expectations on the two sides are on different probability
spaces.) The next theorem from \cite{BKS}
follows easily from (\ref{keyformula}).

\begin{theorem} \label{thmBKS} Given a sequence
$f_n:\{0,1\}^{V_n}\mapsto \{-1,1\}$, the following three conditions
are equivalent. \\
1. The sequence is NS.\\
2. For every $\epsilon > 0$,
$\lim_{n\to\infty}
\E[(1-\epsilon)^{|{\cal S}_{f_n}|}I_{{\cal S}_{f_n}\neq\emptyset}]=0$. \\
3. For every $k$,
\[
\lim_{n\to\infty} \sum_{0<|S|<k}\hat{f_n}(S)^2=0.
\]
(In words, on $\{{\cal S}_{f_n}\neq \emptyset\}$,
$|{\cal S}_{f_n}|\to \infty$ in distribution.)
\end{theorem}

The next theorem, also from \cite{BKS}, similarly
follows easily from (\ref{keyformula}).
\begin{theorem} \label{thmBKSStab} A sequence $f_n:\{0,1\}^{V_n}\mapsto \{-1,1\}$
is NStable iff
\[
\lim_{k \rightarrow \infty}\sup_n \sum_{|S|\geq k} \hat{f}_n(S)^2=0.
\]
(The latter condition means that the family $\{|{\cal S}_{f_n}|\}$ is tight.)
\end{theorem}

\section[Elementary characterizations of Exclusion Sensitivity and Exclusion Stability]{Elementary characterizations of Exclusion Sensitivity and Exclusion Stability}
\label{secXSusingFourier}
The main purpose of this section is to establish elementary necessary and
sufficient conditions for XS and XStable in terms of the Fourier coefficients;
this will also immediately yield
Proposition \ref{propalltsame}. We will also relate a certain strengthening
of these conditions to ordinary noise sensitivity; this will be done in
Proposition \ref{prop:XCgivesNSalmost}.

For the exclusion process
$\{\eta^n_t\}_{t\geq  0}$ with respect to $\{G_n,\alpha_n\}_{n \in \N}$,
it is immediate to check that for $S,S'\subseteq V_n$,
$$
{\mathbb P}(\chi_S(\eta^n_0)=\chi_{S'}(\eta^n_{t}))
=P^n_t(S,S')+1/2(1-P^n_t(S,S'))
$$
which implies that
\begin{equation} \label{eqn10}
{\mathbb E}[\chi_S(\eta^n_0)\chi_{S'}(\eta^n_{t})]
=2{\mathbb P}(\chi_S(\eta^n_0)=\chi_{S'}(\eta^n_{t}))-1
=P^n_t(S,S').
\end{equation}

It is not hard to see that for the SNPS of equation (\ref{eqnNSequiv}), the
functions $\{\chi_S\}$ are also eigenfunctions of the
corresponding Markov semigroup. That is, if $\{\omega_t\}_{t\geq 0}$ is as
in equation (\ref{eqnNSequiv}), we have that
$\E[\chi_S(\omega_t)|\,  \omega_0]=\chi_S(\omega_0)e^{-|S|t}$.
These functions are however not eigenfunctions for the
exclusion process as is easy to see.

The set of eigenfunctions of the Markov semigroup corresponding to the
symmetric exclusion process is in general much harder to describe
(although they are useful, see \cite{Diac}). It will
however be worthwhile to say a few words.
Let $\{\eta^n_t\}_{t\geq 0}$ be the exclusion process on
$\{G_n,\alpha_n\}_{n \in \N}$. We have that
$\E[\chi_S(\eta^n_t)|\,  \eta^n_0]=\sum_{S':|S'|=|S|} P^n_t(S',S)
\chi_{S'}(\eta^n_0)$.
Therefore it is easy to see that for any function $g_k$ which is in
$\Span(\{\chi_S\}_{|S|=k})$,
$\E[g_k(\eta^n_t)|\,  \eta^n_0]$ must also be in $\Span(\{\chi_S\}_{|S|=k})$. Therefore,
there exists a basis of $\Span(\{\chi_S\}_{|S|=k})$ consisting of eigenfunctions of
the restriction of the Markov semigroup to $\Span(\{\chi_S\}_{|S|=k})$.
It is easy to see that $\sum_{|S|=k}\chi_S$ is one such eigenfunction
with eigenvalue 1 since, in this case, $\sum_{|S|=k}\chi_S$ will only be a
function of the number of $1$'s in $\eta^n_t$ which is constant in $t$.

Therefore we can conclude that, for each $n$,
there exists an orthonormal basis
$\{\phi_l^k: 0\leq k \leq |V_n|, 1\leq l \leq {|V_n| \choose k}\}$
of the space of all functions $f:\{0,1\}^{V_n} \mapsto \R$
which are all eigenfunctions for the Markov semigroup corresponding to the
exclusion process with respect to $\{G_n,\alpha_n\}_{n \in \N}$. Furthermore,
for all such $k$ and $l$,
\[
\phi_l^k=\sum_{|S|=k} a_{l,S}\chi_{S},
\]
for some choice of $\{a_{l,S}\}_{|S|=k}$. We have that
\[
\E[\phi_l^k(\eta^n_0)\phi_l^k(\eta^n_t)]=e^{-\lambda_{l}^k t},
\]
for some $\lambda_{l}^k\ge 0$. (The eigenfunctions and eigenvalues depend on
$n$ but we ignore this in the notation.)
Furthermore we assume that for each $k$,
$\{\phi_l^k\}_l$ are ordered so that
$\lambda_{l}^k \le \lambda_{l+1}^k$ for each $l$.
For a function $f:\{0,1\}^{V_n} \mapsto \R$
we let $\hat{f}(\phi_l^k):=\E[f\phi_l^k]$.

For any $C<\infty$, let  $\Lambda_{C}:=\{\lambda_{l}^k:\lambda_{l}^k\leq C\}$
and let $\Phi_C:=\{\phi_l^k: \lambda_{l}^k \leq C\}\setminus \phi_1^0$.
The reason why we disregard $\phi_1^0$ is that this eigenfunction is simply
$\chi_{\emptyset}$, the constant function.

\begin{proposition} \label{prop2}
Let $\{G_n,\alpha_n\}_{n \in \N}$ and $f_n:\{0,1\}^{V_n} \rightarrow \{-1,1\}$
be given. Then we have that
for fixed $t$, the following three conditions are equivalent. \\
1.
\[
\lim_{n\to\infty}
{\mathbb E}[f_n(\eta^n_0)f_n(\eta^n_{t})]-{\mathbb E}[f_n(\eta^n_0)]^2=0.
\]
2.
\[
\lim_{n\to\infty}
\sum_{S\neq \emptyset}\hat{f_n}(S) \sum_{S':|S|=|S'|}\hat{f_n}(S')P^n_{t}(S,S')=0.
\]
3. For any $C<\infty$
\begin{equation} \label{eqn13}
\sum_{\phi_l^k\in\Phi_C} \hat{f_n}(\phi_l^k)^2 \rightarrow 0.
\end{equation}
\end{proposition}

\noindent
{\bf Remark:}
The equivalence of 1 and 3 here is analogous to the
equivalence between $1$ and $3$ in Theorem \ref{thmBKS}.\\

\noindent
{\bf Proof.}
We first prove the equivalence of 1 and 2.
We have that
\begin{eqnarray}
\lefteqn{{\mathbb E}[f_n(\eta^n_0)f_n(\eta^n_{t})]} \label{eqn2} \\
& & ={\mathbb E}[\sum_{S} \hat{f_n}(S)
\chi_S(\eta^n_0)\sum_{S'} \hat{f_n}(S')\chi_{S'}(\eta^n_{t})] \nonumber \\
& & =\sum_{S}\hat{f_n}(S) \sum_{S':|S|=|S'|}\hat{f_n}(S'){\mathbb E}[\chi_S(\eta^n_0)\chi_{S'}(\eta^n_{t})] \nonumber \\
& & =\sum_{S}\hat{f_n}(S) \sum_{S':|S|=|S'|} \hat{f_n}(S')P^n_{t}(S,S').\nonumber
\end{eqnarray}
Of course we have that
\[
\E[f_n]=\hat{f_n}(\emptyset),
\]
so that
\[
{\mathbb E}[f_n(\eta^n_0)f_n(\eta^n_{t})]-{\mathbb E}[f_n(\eta^n_0)]^2
=\sum_{S\neq \emptyset}\hat{f_n}(S) \sum_{S':|S|=|S'|}\hat{f_n}(S')P^n_{t}(S,S').
\]
Hence 1 and 2 are equivalent.

We now prove the equivalence of 1 and 3. Expanding $f_n$, we get that
\[
f_n=\sum_{k=0}^{|V_n|} \sum_{l=1}^{{|V_n| \choose k}} \hat{f_n}(\phi_l^k)\phi_l^k.
\]
Therefore
\begin{eqnarray*}
\lefteqn{\E[f_n(\eta^n_0) f_n(\eta^n_t)]-\E[f_n(\eta^n_0)]^2}\\
& & =\sum_{k=1}^{|V_n|} \sum_{l=1}^{{|V_n| \choose k}} \hat{f_n}(\phi_l^k)^2 \E[\phi_l^k(\eta^n_0) \phi_l^k(\eta^n_t)] \\
& & =\sum_{k=1}^{|V_n|} \sum_{l=1}^{{|V_n| \choose k}} \hat{f_n}(\phi_l^k)^2 e^{-\lambda_{l}^k t}
\geq e^{-Ct}\sum_{\phi_l^k\in\Phi_C} \hat{f_n}(\phi_l^k)^2,
\end{eqnarray*}
so that $\{f_n\}_{n\in\N}$ cannot satisfy 1 if equation (\ref{eqn13})
fails for some $C<\infty$. Hence 1 implies 3. Similarly,
\[
\E[f_n(\eta^n_0) f_n(\eta^n_t)]-\E[f_n(\eta^n_0)]^2 \leq \sum_{\phi_l^k\in\Phi_C} \hat{f_n}(\phi_l^k)^2+
\sum_{\phi_l^k\in (\Phi_C)^c \setminus \{\phi_1^0\}} \hat{f_n}(\phi_l^k)^2 e^{-C t}.
\]
Then the RHS can be made arbitrarily small for $n$ sufficiently large,
if equation (\ref{eqn13}) holds for every $C<\infty$. Hence
3 implies 1.
\fbox{}\\

\medskip\noindent
{\bf Proof of Proposition \ref{propalltsame}}
The proof is immediate from Proposition \ref{prop2} since
condition 3 does not depend on $t$.
\fbox{}\\

There is an elementary characterization of NStable as well.
Since the proofs follow the same lines, we leave them to the reader.

\begin{proposition} \label{prop2stable}
Let $\{G_n,\alpha_n\}_{n \in \N}$ and $f_n:\{0,1\}^{V_n} \rightarrow \{-1,1\}$
be given. Then the following three conditions are equivalent. \\
1. $\{f_n\}_{n \in \N}$ is XStable.  \\
2. for any $\delta>0$, there exists a $t>0$ such that for every $n$,
\[
\sum_{S}\hat{f_n}(S) \sum_{S':|S|=|S'|}\hat{f_n}(S')P^n_{t}(S,S')\geq 1-\delta.
\]
3. For any $\delta>0$, there exists $C$ such that for every $n$,
\begin{equation} 
\hat{f_n}(\phi_1^0)^2 +
\sum_{\phi_l^k\in\Phi_C} \hat{f_n}(\phi_l^k)^2 \ge 1-\delta.
\end{equation}
\end{proposition}

\noindent
{\bf Remark:}\\
The equivalence of 1 and 3 is analogous to Theorem \ref{thmBKSStab}.

\bigskip

At this point it is very natural to ask the following:
does there exist a sequence $\{G_n,\alpha_n\}_{n \in \N}$ and a sequence
$f_n:\{0,1\}^{V_n}\mapsto \{-1,1\}$ for which
\begin{equation} \label{eqn1}
\sum_{S\neq \emptyset}\hat{f_n}(S) \sum_{S':|S|=|S'|}\hat{f_n}(S')P^n_{t}(S,S')
\rightarrow 0
\end{equation}
but
\begin{equation} \label{eqn4}
\sum_{S\neq \emptyset}|\hat{f_n}(S)| \sum_{S':|S|=|S'|}|\hat{f_n}(S')|P^n_{t}(S,S')
\rightarrow 0 \ \
\end{equation}
fails?  The reason it is natural to ask this question is that one
would want to know if XS sometimes occurs due to cancellation of
the positive and negative terms in (\ref{eqn1}). We will see below
that it can in fact happen that (\ref{eqn1}) holds but (\ref{eqn4}) fails,
showing that cancellation of terms is relevant.
However, before doing that, we point out the interesting fact that had it been
the case that (\ref{eqn1}) implied (\ref{eqn4}), this would have
provided an alternative proof to part 1 of Theorem \ref{thm2}.
This is stated in the
following easy proposition which might be interesting in its own right.

\begin{proposition} \label{prop:XCgivesNSalmost}
Given a sequence of Boolean functions, if (\ref{eqn4}) holds, then the
sequence is NS.
\end{proposition}
\noindent
{\bf Proof.}
We first observe that with $\epsilon=1-e^{-t}$
\begin{eqnarray*}
\lefteqn{\sum_{S\neq \emptyset}|\hat{f_n}(S)| \sum_{S':|S|=|S'|}|\hat{f_n}(S')|P^n_{t}(S,S')}\\
& & \geq \sum_{S\neq \emptyset}\hat{f_n}(S)^2P^n_t(S,S)
\geq \sum_{S\neq \emptyset}\hat{f_n}(S)^2 e^{-t|S|}
=\E[(1-\epsilon)^{|{\cal S}_{f_n}|}I_{{\cal S}_{f_n}\neq\emptyset}]
\end{eqnarray*}
where the last inequality follows from the assumption on the rates $\alpha_n$.
Now apply Theorem \ref{thmBKS}.
\fbox{}\\

Concerning the question of whether (\ref{eqn4}) can
hold for the complete graph dynamics, it is easy to check
that Example \ref{ex7} yields such an example.

\medskip

Example \ref{ex2} below shows that (\ref{eqn1}) does not imply
(\ref{eqn4}) in general. Before proceeding with this example,
we need to make an observation. For $B\subseteq [n]$, let
$\sigma_B:\{0,1\}^{V_n}\rightarrow \{0,1\}^{V_n}$ be such
that for $\omega \in \{0,1\}^{V_n}$
\[
\sigma_B(\omega)_i=\left\{\begin{array}{ccc}
(\omega_i+1)  \textrm{ mod }\, 2,& \textrm{if} & i\in B\\
\omega_i & \textrm{otherwise.} \\
\end{array} \right.
\]
In words, $\sigma_B$ just flips the bits in $B$.

It is easy to see that
\begin{equation} \label{switchformula}
\chi_S(\sigma_B(\omega))=(-1)^{|S\cap B|}\chi_S(\omega).
\end{equation}

\begin{example} \label{ex2}
\end{example}
Here $|\omega|$ refers to the number of 1's in $\omega$.
Let $G_n=(V_n,E_n)$ be the graph consisting of $n$ isolated edges and
let $\alpha_n=1$ for every $n$. Let $\omega \in \{0,1\}^{V_n}$
have distribution $\pi_{1/2}$.  Define the Boolean function $f_n$
in the following way. We let
\[
f_n(\omega)=\left\{ \begin{array}{cc}
1 & \textrm{ if } |\omega|\in\{4k,4k+1\}, \ \mbox{ for some } k \in
\{0,1,\ldots, n/2\} \\
-1 & \textrm{ otherwise.}
\end{array}
\right.
\]
Arbitrarily, order the edges of $G_n$, $e_1,e_2,\ldots,e_n$ and denote
the two endpoints of $e_k$ by $v_{k,1}$ and $v_{k,2}$.
Let $B_n:=\{v_{k,2}: 1\leq k \leq n\}$ and define
$g_n:=f_n \circ \sigma_{B_n}$.

Using (\ref{switchformula}), we have
\[
g_n(\omega)=f_n(\sigma_{B_n}(\omega))= \sum_S \hat{f_n}(S)\chi_S(\sigma_{B_n}(\omega))=
\sum_S (-1)^{|S\cap B_n|}\hat{f_n}(S)\chi_S(\omega),
\]
so that $\hat{g_n}(S)=(-1)^{|S\cap B_n|}\hat{f_n}(S)$.

Trivially $\{f_n\}_{n\in \N}$ will not be XS due to particle conservation.
However, $g_n$ will be XS as can be seen as follows. If
$\eta_0(v_{k,1})=\eta_0(v_{k,2})$, then $\sigma_{B_n}(\eta_t(v_{k,i}))$
is constant in time for $i=1,2$. If instead
$\eta_0(v_{k,1}) \neq \eta_0(v_{k,2})$,
then $\sigma_{B_n}(\eta_0)(v_{k,1})=
\sigma_{B_n}(\eta_0)(v_{k,2})$ and will be
$1$ or $0$ with equal probability. Furthermore,
$(\sigma_{B_n}(\eta_t)v_{k,1},\sigma_{B_n}(\eta_t)v_{k,2})$
will change between the states $\{(0,0), (1,1)\}$ with rate 1.
Using the fact that the number of $k$ such that
$\eta_0(v_{k,1}) \neq \eta_0(v_{k,2})$ has a binomial distribution
with parameters $n$ and $1/2$, one can easily show
that $\{g_n\}_{n\in\N}$ is XS. The details are left to
the reader.

Since $f_n(\omega)$ is constant under permutations of $\omega$, it
is easy to see that $\hat{f_n}(S)=\hat{f_n}(S')$ for every
$|S|=|S'|$. Therefore,
\[
|\hat{g_n}(S)\hat{g_n}(S')|=|\hat{f_n}(S)\hat{f_n}(S')|=\hat{f_n}(S)\hat{f_n}(S'),
\]
for every $S,S'$ with $|S|=|S'|$.
Since $\{g_n\}_{n\in \N}$ is XS, Proposition \ref{prop2} tells us that
(\ref{eqn1}) must hold for this sequence. However,
(\ref{eqn4}) cannot hold for this sequence since if it did, then by
the above, (\ref{eqn1}) would hold for $\{f_n\}_{n\in \N}$,
contradicting the fact that the latter is not XS (again using
Proposition \ref{prop2}).
\vspace{3mm}

\section{Some further examples}
\label{prop:CGXSNOTIMPLY2-NNXSSECTION}

Given an exclusion process on $V$ and $S\subseteq V$, recall that $S_t=\pi_t(S)$.
It will be useful to know something about the distribution of the size of
$S_t\cap S$. If we consider the complete graph exclusion process on
$V$, it is easy to see that conditioned on $|S_t\cap S|$, the
set $S \setminus S_t$ will be uniformly distributed on $S$ and likewise
$S_t\setminus S$ will be uniformly distributed on $S^c:=V\setminus S$.
We will need to analyse the distribution of $|S_t\cap S|$ for $|S|$ small.
This will be done in Lemma \ref{lemma3} below, which we will also use in the
proof of Theorem \ref{thm5} later on.

For integer valued random variables $X$ and $Y$, we will write
$X \preceq Y$ if we can couple $X$ and $Y$ so that $\P(X \leq Y)=1$.

\begin{lemma} \label{lemma3}
If $|S|<|V|/2$ then for all $t\ge 0$,
${\rm Bin}(|S|,\epsilon) \preceq |S_t\setminus S|$, where
\[
\epsilon=\left(1-e^{-\left(1-|S|/|V|\right)t}\right)\left(1-\frac{|S|}{|V|-|S|}\right).
\]
\end{lemma}
\noindent
{\bf Proof.}
It will be convenient to reformulate the dynamics of the complete
graph exclusion process.
Assume that $\eta_0(S)\equiv 1$ and that $\eta_0(S^c)\equiv 0$;
we will think of
$S$ initially being occupied by particles while $S^c$ is empty.
Let $x_1,\ldots,x_{|S|}$ be the particles of $S$. Associate an independent
Poisson process with rate $1-|S|/|V|$ to each of these particles and
when the Poisson clock associated to a particle $x_i$ rings,
the particle chooses uniformly at random a site in $V$ that is unoccupied,
and moves there. It is easy to see that this is equivalent
to the complete graph exclusion dynamics.

Fix $t\ge 0$ and $S$ with $|S|<|V|/2$.
To each particle $x_i$ associate an independent
sequence $\{U_{i,j}\}_{j\geq 1}$ of independent
$U[0,1]$ random variables. Assume that the Poisson clock associated
to $x_i$ rings at some time $\tau_i \in [0,t]$ and let $j(\tau_i)$ be
the smallest $j$ such that $U_{i,j}$ has not been used previously.
If $U_{i, j(\tau_i)}\leq 1-|S|/(|V|-|S|)$, choose an unoccupied
site of $S^c$ uniformly at random and move $x_i$ there. If instead $U_{i, j(\tau_i)}> 1-|S|/(|V|-|S|)$,
then move $x_i$ to an unoccupied site (uniformly) in $S$ with probability
$|S\setminus S_{\tau_i-}|/|S|$ and otherwise jump to an unoccupied site (uniformly) in $S^c$.
Obviously, the probability that $x_i$ jumps to some unoccupied site in $S$ is
\[
\P(U_{i, j(\tau_i)}> 1-|S|/(|V|-|S|))\frac{|S\setminus S_{\tau_i-}|}{|S|}=\frac{|S\setminus S_{\tau_i-}|}{|V|-|S|},
\]
and so the use of $U_{i,j(\tau_i)}$ correctly describes the dynamics.

Let $J(i)$ be the largest $j$ such that $U_{i,j}$ was used up to time $t$
(for particle $x_i$).
We say that a particle is ``good'' if the Poisson clock associated to $x_i$
rings at least once and $U_{i,J(i)}\leq 1-|S|/(|V|-|S|)$. Clearly the
$x_i$ will be good independently of each other. Furthermore,
\[
\P(x_i \textrm{ is good})=\left(1-e^{-\left(1-|S|/|V|\right)t}\right)\left(1-\frac{|S|}{|V|-|S|}\right).
\]
Finally it is easy to see that if $x_i$ is good then
$x_i\in S_t\setminus S$.
\fbox{}\\

\noindent
{\bf Remark:} It might seem that this argument was slightly more
complicated than is needed. For example, one might think that it is
easy to show that conditioned on the event that the points in $S$
which have been updated by time $t$ is exactly some subset $S'$, then
the conditional distribution of $S_t$ is the union of $S\setminus S'$
and a uniform subset of $S^c\cup S'$ of size $|S'|$. However,
interestingly, this is false. This is why a slightly more involved
argument was needed.

\medskip\noindent
Before we proceed, we make the following observations.
If $\omega,\omega^S$ are $\pi_{1/2}$-distributed and such that they agree on $S^c$ and are independent on
$S$, then it is not hard to check that
\[
\E[f(\omega)f(\omega^S)]=\sum_{S'\cap S = \emptyset }\hat{f}(S')^2.
\]
It follows that for any $T_1 \subseteq T_2$,
\begin{equation} \label{eqn9}
\E[f(\omega)f(\omega^{T_1})]\geq \E[f(\omega)f(\omega^{T_2})].
\end{equation}
Letting $g(S):=\E[f(\omega)f(\omega^S)]$ and noting that by (\ref{eqn9}),
 $g(S)$ is decreasing in $S$, we conclude that if $T_1,T_2$ are random
sets such that $T_1 \subseteq T_2$ with
probability one, then (\ref{eqn9}) also holds for $T_1,T_2$.
We implicitly assume that the randomness determining the random sets
$T_1$ and $T_2$ is independent of $\omega$ and any rerandomization of
$\omega$.

We now give our first example of this section which is the following.
\begin{example} \label{exCGXSnotNN}
Consider the boxes $\Lambda^2_n:=\{1,\ldots,n\}^2$ and $\Lambda^2_{n^2}:=\{1,\ldots,n^2\}^2$.
Assume for simplicity that $n$ is odd and define for $k\in \Z$
$I_k:=[\frac{n^2}{2}+k n^{2/3},\frac{n^2}{2}+(k+1) n^{2/3} ]$. Let for $\omega\in\{0,1\}^{\Lambda^2_{n^2}}$
\[
f_n(\omega)=\left\{
\begin{array}{cc}
1 & \textrm{ if } |\omega_{\Lambda^2_n}|\in  I_k, \ \ k \textrm{ odd}\\
-1 & \textrm{ otherwise, }
\end{array}\right.
\]
where $\omega_{\Lambda^2_n}$ denotes the restriction of $\omega$ to $\Lambda^2_n$. By symmetry,
using the fact that $n$ is odd, it is easy to see that $\E[f_n(\omega)]=0$.
\end{example}

\begin{proposition} \label{prop:CGXSNOTIMPLY2-NNXS}
Example \ref{exCGXSnotNN} is CGXS but not 2-NNXS.
In fact, this sequence is 2-NNXStable.
(The example is easily generalized to yield, for any $d$, a sequence of Boolean functions on
$\Lambda_n^d$ which is CGXS and d-NNXStable.)
\end{proposition}

\noindent
{\bf Proof.}
We will start by showing that $f_n$ is not 2-NNXS. Let $S=\Lambda^2_n$ and $S_t=\pi^{n,NN}_t(S)$,
where the superscript $NN$ stresses that we are considering the nearest neighbor model.
It is easy to see that
$|S \setminus S_t |$ is bounded from above by the total number of times that the Poisson
clocks associated to edges between
$\Lambda^2_n$ and $\Lambda^2_{n^2} \setminus \Lambda^2_n$ ring. This is trivially a Poisson distributed number of times with parameter
$n t /2$. Let $|\eta^n_0(\Lambda^2_n)|$ and $|\eta^n_t(\Lambda^2_n)|$ denote the number
of ones in $\Lambda^2_n$ at time 0 and $t$ respectively.
Let $\delta>0$, by the above observation it is easy to see that there exists a constant $0<C<\infty$ such that
$\P(||\eta^n_0(\Lambda^2_n)| - |\eta^n_t(\Lambda^2_n)| | \geq C \sqrt{n}) \leq \delta$
for all $n$ large and $t\in [0,1]$.
Furthermore, it is straightforward to check that
\[
\P\left(|\eta^n_0(\Lambda^2_n)|\in \left[\frac{n^2}{2}+k n^{2/3}-C\sqrt{n},\frac{n^2}{2}+k n^{2/3}+C\sqrt{n}\right] \textrm{ for some } k\in \Z\right)
\leq \delta
\]
for $n$ large enough. Therefore, $\P(f_n(\eta^n_0)=f_n(\eta^n_t)) \geq 1-2 \delta$
and so $f_n$ is not 2-NNXS and in fact is now easily seen to be 2-NNXStable.

We proceed to show that $f_n$ is CGXS.
Let $S=\Lambda^2_n$ and $S_t=\pi^{n,CG}_t(S)$, where $CG$ refers to the complete graph
exclusion dynamics.
According to Lemma \ref{lemma3}, when running the
complete graph dynamics, $|S \setminus S_t |$ dominates a
${\rm Bin}(n^2,\epsilon)$ random variable with
\[
\epsilon=\left(1-e^{-\left(1-n^2/n^4\right)t}\right)\left(1-\frac{n^2}{n^4-n^2}\right)
>1-e^{-t/2},
\]
where the inequality holds for all $n$ large enough. Conditioned on $|S \setminus S_t |$
the distribution of the set $T_2:=S \setminus S_t$ is uniform among
all subsets of $S$ that are of size $|S \setminus S_t |$.
Letting $T_1$ be such that $x\in T_1$ for every $x\in S$ independently with probability $\epsilon$,
we can couple the two random variables $T_1$ and $T_2$
such that $T_1 \subseteq T_2$. Therefore, for our very particular choice of
function $f_n$, which is invariant under permutations leaving $\Lambda^2_n$ fixed,
\begin{eqnarray*}
\lefteqn{N(f_n,t,G_n,\alpha_n)=\E[f_{n}(\omega)f_{n}(\omega^{T_2})]-\E[f_{n}(\omega)]^2}\\
& & \leq \E[f_{n}(\omega)f_{n}(\omega^{T_1})]-\E[f_{n}(\omega)]^2
= \E[f_{n}(\omega)f_{n}(\omega^{\epsilon})]-\E[f_{n}(\omega)]^2,
\end{eqnarray*}
according to (\ref{eqn9}) and the discussion thereafter. Hence, if $f_n$ is NS it must be CGXS.

Showing that $f_n$ is NS is a straightforward exercise. Observe that we resample $N \sim {\rm Bin}(n^2,\epsilon)$
bits and therefore with high probability $|\omega|-|\omega^{\epsilon}|$ will be of order $O(n)$.
Since $f_n$ changes value in intervals of length $n^{2/3}$,
NS follows quite easily.
\fbox{}\\

\noindent
{\bf Remark:}\\
Observe that we not only showed 2-NNXStability but showed the stronger fact that
\[
\lim_{n\rightarrow \infty} \max_{t\in [0,1]}\P[f_n(\eta^n_0)\neq f_n(\eta^n_{t})]=0.
\]
Such behavior cannot occur when studying usual noise sensitivity except in
degenerate situations.

The following example demonstrates that this same behavior can also occur
in the monotone case. In fact, in the example below, we show NS rather than
CGXS but an appeal to Corollary \ref{cor:monotone} allows us to conclude
CGXS.

\medskip\noindent
\begin{example} \label{ex5} (Crossings on majority)
\noindent Consider $\Lambda^2_{n} \subseteq \Z^2$. Partition this box into
subboxes of sidelength $n^{\alpha}$ (where $\alpha>0$ and $n^{\alpha}$ is
assumed for simplicity to be an odd integer).
In every one of these subboxes $B$, let $g_B$ be
the majority function on that subbox. Consider now the event
${\cal C}_n$ that there exists a path from left to right of
$\Lambda_{n}^2$ by subboxes with $g_B=1$.
(The probability of such an event is the same as the probability
of having an open crossing of an $n^{1-\alpha}\times n^{1-\alpha}$
box.)
\end{example}

\medskip\noindent
\begin{proposition} \label{acounterexample}
For all $\alpha \in (0,1)$, Example \ref{ex5} is NS and for $\alpha > 8/9$,
it is not 2-NNXS and is in fact 2-NNXStable.
\end{proposition}

\noindent
{\bf Proof of Proposition \ref{acounterexample}.} \\
We first show NS.
Fix $\epsilon > 0$. It is well known that if we stick to functions
with mean $0$, the correlation between $f(\omega)$ and
$f(\omega^\epsilon)$ is maximized when $f$ is a so-called dictator
function which means that it only depends on one bit in which
case the correlation is $1-\epsilon$. Hence after the rerandomization,
we can view the state of each subbox as having been
rerandomized with a probability of at least $\epsilon$. Since it
is known that crossings of percolation are noise sensitive
(proved in \cite{BKS}), it follows that this example is NS as well.

We now show that 2-NNXS fails for $\alpha>8/9$. Fix $t=1$. Fix a subbox $B$ and
let $X$ be the number of particles outside $B$ at time 0 which are in $B$
at time 1. Obviously $X$ is dominated by a Poisson distributed random variable $Y$
with parameter $n^{\alpha}$. Therefore
$\P(X \geq n^{3\alpha/2})\leq \P(Y \geq n^{3\alpha/2}) \leq n^{-\alpha/2}$ by Markov's inequality.
By (\ref{eqn9}) and the discussion thereafter, it follows as in the proof of
Proposition \ref{prop:CGXSNOTIMPLY2-NNXS} that
\[
{\mathbb P}(g_B(\eta^n_1)\neq g_B(\eta^n_0)) \leq {\mathbb P}(g_B(\eta^n_1)\neq g_B(\eta^n_0)|X=n^{3\alpha /2})
+n^{-\alpha/2}.
\]
If we were to rerandomize the bits independently with probability $\epsilon_n=n^{-\alpha/2}$,
the number $Z$ of rerandomized bits would be ${\rm Bin}(n^{2\alpha}, n^{-\alpha/2})$ distributed.
As above, using (\ref{eqn9}) and the discussion thereafter in the first inequality
\begin{eqnarray*}
\lefteqn{\P(g_B(\eta^n_1)\neq g_B(\eta^n_0)|X=n^{3\alpha /2})}\\
& & =\P(g_B(\omega)\neq g_B(\omega^{\epsilon_n}) |\,  Z = n^{3\alpha /2}) \\
& & \leq \P(g_B(\omega)\neq g_B(\omega^{\epsilon_n}) |\,  Z \geq n^{3\alpha /2})\\
& & \leq \frac{\P(g_B(\omega)\neq g_B(\omega^{\epsilon_n}))}{\P(Z \geq n^{3\alpha /2})}
=O(1)\P(g_B(\omega)\neq g_B(\omega^{\epsilon_n})).
\end{eqnarray*}
In \cite{BKS}, Remark 3.6 they conclude that $\P(g_B(\omega)\neq g_B(\omega^{\epsilon_n}))\leq O(1)\epsilon_n^{1/2}$.
We therefore conclude that
\[
{\mathbb P}(g_B(\eta^n_1)\neq g_B(\eta^n_0))\leq O(1)n^{-\alpha/4}.
\]
It follows that the probability that there is some subbox whose
$g_B$ value changes is at most
$$
O(1)n^{2-2\alpha}/n^{\alpha/4}.
$$
If $\alpha > 8/9$, this approaches 0. This implies we cannot be 2-NNXS and it is
also easily argued that we have  2-NNXStability.
\fbox{}\\

\medskip\noindent
{\bf Remark:}  We again have a strengthening of stability which was
described in the remark following the proof of Proposition
\ref{prop:CGXSNOTIMPLY2-NNXS}.

\section[Exclusion Sensitivity implies Noise Sensitivity]{Exclusion Sensitivity (XS) implies Noise Sensitivity (NS)} \label{secproofs}

\medskip\noindent
{\bf Proof of Theorem \ref{thm2}.}
We start by proving the first statement.
Assume that the sequence $\{f_n\}_{n\in\N}$ is not NS.
From Theorem \ref{thmBKS} it follows that there exists a $k<\infty$
and $c>0$ such that for some subsequence $\{n_l\}_{l\in\N}$
\[
\sum_{0<|S|\leq k} \hat{f}_{n_l}(S)^2\geq c.
\]
For notational simplicity we will assume that in fact this holds for
every $n$.
Next, let $t>0$ be such that $P^n_{t}(S,S) \geq 3/4$ for every $S$ with
$|S|\leq k$ and
for all $n$. This can be done by choosing
$t$ such that $e^{-kt}\geq 3/4$, by our assumptions on $\alpha_n$.

Let $P^{n,k}_{t}$ denote the restriction of $P^n_t$ to the nonempty
sets of size at most $k$. Note that due to particle conservation,
this is a transition matrix. Since the diagonal elements are all at least
$3/4$, there is a transition matrix $\tilde{P}^{n,k}(S,S')$ for
$0<|S|,|S'|\le k$, so that
\[
P^{n,k}_{t}(S,S')=\frac{3}{4}I+\frac{1}{4}\tilde{P}^{n,k}(S,S'),
\]
where $I$ is the appropriate identity matrix.
Clearly, all of the matrices involved are also symmetric matrices.
Therefore, as all of the eigenvalues of $\tilde{P}^{n,k}$
are real and sit inside $[-1,1]$, we easily get that the minimum
eigenvalue of $P^{n,k}_{t}$ is at least $1/2$.

Let $\hat{f_n}^k$ be a vector consisting of $\hat{f_n}(S)$ for every
$0<|S|\leq k$. Since $P^{n,k}_{t}$ is symmetric, we can write
$\hat{f_n}^k=\sum \alpha_i v_i$ where $\{v_i\}$ is
an orthonormal basis of $P^{n,k}_{t}$.
Letting $\lambda_i$ denote the eigenvalue of the eigenvector $v_i$, we get that
\begin{eqnarray} \label{eqn5}
\lefteqn{\sum_{0<|S|\leq k}\hat{f_n}(S) \sum_{S':|S|=|S'|}\hat{f_n}(S')P^{n,k}_{t}(S,S')}\\
& & = \sum_{i} \lambda_i \alpha_i^2 \geq \frac{1}{2}\sum_{i} \alpha_i^2
=\frac{1}{2}\sum_{0<|S|\leq k} \hat{f_n}(S)^2 \geq \frac{c}{2}. \nonumber
\end{eqnarray}

On the other hand, since all eigenvalues of $P^n_{t}(S,S')$ are nonnegative, we also have that
\begin{equation}
\sum_{|S|>k}\hat{f_n}(S) \sum_{S':|S|=|S'|}\hat{f_n}(S')P^n_{t}(S,S')\geq 0,
\nonumber
\end{equation}
giving that
\begin{equation}
\sum_{|S|>0}\hat{f_n}(S) \sum_{S':|S|=|S'|}\hat{f_n}(S')P^n_{t}(S,S')
\ge c/2.
\nonumber
\end{equation}
Hence, $\{f_n\}_{n\in \N}$ is not XS by Proposition \ref{prop2}.

For the second statement, assume that the sequence is NStable. Then, given
$\delta>0$, Theorem \ref{thmBKSStab} yields that
there exists a $k<\infty$ such that for every $n$
\[
\sum_{0\leq |S|<k} \hat{f}_n(S)^2\geq 1-\delta.
\]
Furthermore, by choosing $t>0$ sufficiently small we can (as above)
ensure that the smallest eigenvalue of $P^{n,k}_{t}$ is at least $1-\delta$
uniformly in $n$. As above, we can conclude that
\[
\sum_{S}\hat{f_n}(S) \sum_{S':|S|=|S'|}\hat{f_n}(S')P^n_{t}(S,S')\geq 1-2\delta.
\]
Therefore $\{f_n\}_{n\in \N}$ is XStable by Proposition \ref{prop2stable}.
\fbox{}\\

\medskip

It is also worthwhile to observe that a similar argument to the proof of
Theorem \ref{thm2} shows that if the sequence $\{f_n\}_{n\in\N}$ is XS
with respect to a certain $\{G_n,\alpha_n\}_{n \in \N}$,
then in fact we must have that for any $k<\infty$,
\[
\lim_{n\rightarrow \infty} \sum_{|S|\geq |V_n|-k}\hat{f_n}(S)^2=0.
\]
This is another way to see why $f_n=\chi_{V_n}$ cannot be XS since here $\hat{f_n}(V_n)=1$
for any $n$ and any sequence $\{G_n, \alpha_n\}_{n\in\N}$. Hence we see that a necessary condition for
$\{f_n\}_{n\in \N}$ to be XS (with respect to any $\{G_n, \alpha_n\}_{n\in\N}$)
is that for any $k$
\[
\lim_{n\rightarrow \infty} \sum_{0<|S|\leq k}\hat{f_n}(S)^2+ \sum_{|S|\geq |V_n|-k}\hat{f_n}(S)^2=0.
\]

In fact one might ask whether there exists an example $\{f_n\}_{n\in \N}$ which is not
XS (with respect to any $\{G_n, \alpha_n\}_{n\in\N}$) but for which
\[
\lim_{n\rightarrow \infty} \sum_{0<|S|\leq k}\hat{f_n}(S)^2+ \sum_{|S|\geq |V_n|-k}\hat{f_n}(S)^2=0.
\]
The following provides such an example.
\begin{example} \label{ex4}
\end{example}
Let $|V_n|=n$ and
\[
f_n(\omega):=\left\{\begin{array}{cc}
1 & \textrm{ if } |\omega|\in[2k \log n,(2k+1)\log n ), \textrm{ for some } 0\leq k \leq n/(2 \log n),\\
-1 & \textrm{ otherwise. }
\end{array}
\right.
\]
Trivially $\{f_n\}_{n\in \N}$ is not XS due to particle conservation. Next,
as the reader can easily check, $\{f_n\}_{n\in \N}$ is NS and therefore for
any $k$,
\[
\lim_{n\rightarrow \infty} \sum_{0<|S|\leq k}\hat{f_n}(S)^2=0.
\]
Furthermore observe that
\[
\hat{f}_n(S)=\E[f_n \chi_S]=\E[f_n \chi_{V_n} \chi_{S^c}]=\hat{g}_n(S^c),
\]
where $g_n:=f_n\chi_{V_n}$. The reader can easily check that $g_n$ is
also NS and furthermore that
\[
\lim_{n\rightarrow \infty} \hat{g_n}(\emptyset)=0.
\]
Therefore for any $k$,
\[
\lim_{n\rightarrow \infty} \sum_{|S|\geq |V_n|-k}\hat{f_n}(S)^2=\lim_{n\rightarrow \infty} \sum_{0\leq |S|\leq k}\hat{g_n}(S)^2=0.
\]

\bigskip

\section[$I\!I(f_n)\to 0$ implies Complete Graph Exclusion Sensitivity]{$I\!I(f_n)\to 0$ implies CGXS and CGXStable implies NStable for monotone functions}
\label{sect:monotonicity}

We will start by presenting the key lemma in proving Theorem
\ref{thmmonNSimpliesCGXS} which is also interesting in itself.
We assume throughout this section that $V_n=[n]$ for every $n$.

\begin{lemma}\label{lemma12}
For every $n$, let $(\xi_n,\xi^*_n)$ be a pair of random
configurations in $\{0,1\}^{[n]}$ satisfying the following conditions:

\begin{enumerate}
\item $\xi_n$ is uniformly distributed on  $\{0,1\}^{[n]}$.

\item The joint distribution of $(\xi_n,\xi^*_n)$ is invariant under all
permutations of $[n]$.

\item For every $\delta>0$ there exists a $C$ such that for every $n$,
\[
\P(d(\xi_n,\xi^*_n)\geq C \sqrt{n})<\delta,
\]
(i.e. the sequence $\{d(\xi_n,\xi^*_n)/\sqrt{n}\}_{n \in \N}$ is tight)
where $d(\xi_n,\xi^*_n):=\sum_{j\in [n]}|\xi(j)-\xi^*(j)|$ is the so-called Hamming distance.

\end{enumerate}
For any such sequence of pairs $(\xi_n,\xi^*_n)$,
if the sequence $\{f_n\}_{n \in \N}$ satisfies $\lim_{n \rightarrow \infty} I\!I(f_n)=0$, then
\[
\lim_{n \rightarrow \infty} \P(f_n(\xi_n)\neq f_n(\xi^*_n))=0.
\]
In particular, this holds if the $\{f_n\}_{n \in \N}$ are monotone and noise sensitive.
\end{lemma}

\noindent
{\bf Proof.} Without loss of generality, we assume $n$ is even. Furthermore, we will
write $(\xi,\xi^*)$ for $(\xi_n,\xi^*_n)$. From $(\xi,\xi^*)$ and additional randomness,
we will construct two random variables whose joint distribution is the same as $(\xi,\xi^*)$.
We will see later the usefulness of this construction.
Given $(\xi,\xi^*)$, let $M_{01}$ be the random number of 0's in
$\xi$ that are changed to 1 in $\xi^*$ and let
$M_{10}$ be defined analogously. We now define a random path through
the hypercube which first moves ``up'' (meaning 0's change to 1's)
and then moves ``down'' (meaning 1's change to 0's). Informally, we choose
$(\xi,M_{01},M_{10})$ as above
and then choose a completely uniform random path starting from
$\xi$ moving ``up'' $M_{01}$ steps and then moving ``down'' $M_{10}$ steps
but where we don't allow that any coordinate changes more than once.
More formally, first choose $(\xi,M_{01},M_{10})$  as above.
Then define a random path $\xi_0,\xi_1,\ldots,\xi_{M_{01}+M_{10}}$ in the hypercube
in the following way.
Let $\xi_0=\xi$, and then for every $\xi_i$ with
$1\leq i \leq M_{01}-1$, $\xi_{i+1}$ is obtained by
picking one of the 0's in $\xi_i$ uniformly at random and then
changing it to a 1. This defines the first part of the
path $\xi_0,\xi_1,\ldots,\xi_{M_{01}}$.
Next, for every $\xi_i$ with $i\in [M_{01},M_{01}+M_{10}-1]$,
$\xi_{i+1}$ is obtained by picking one of the 1's in $\xi_i$,
which was also a 1 for $\xi_{0}$, uniformly at random and then
changing it to a 0. Clearly the distribution of this
random path of random length in the hypercube is invariant under all
permutations. Call this random path $\pi$. Observe crucially that,
due to condition 2, the joint
distributions of $(\xi,\xi^*)$ and $(\xi_0,\xi_{M_{01}+M_{10}})$ are equal.

The \emph{edge boundary} of a subset $A$ of the hypercube $\{0,1\}^{[n]}$
is the set of edges with exactly one endpoint in $A$.
Let $E_n$ be the edge boundary of $\{\xi: f_n(\xi)=1\}$. We have
$$
\P(f_n(\xi)\neq f_n(\xi^*))= \P(f_n(\xi_0)\neq f_n(\xi_{M_{01}+M_{10}}))
\le \P(E_n\cap \pi \neq \emptyset)
$$
where $\pi$ is identified with its set of edges.

For every $k\in \Z$, let $S_k:=\{\xi\in \{0,1\}^{[n]}: |\xi|=n/2+k\}$.
Observe that the number of edges from $S_k$ to $S_{k+1}$ is exactly
${n \choose n/2+k}(n/2-k)$. It follows, from invariance and the fact that
$\pi$ first goes up and then goes down, that
for a fixed edge $e$ in the hypercube between $S_k$ and $S_{k+1}$, we have
$\P(e\in \pi)$ is at most $2\left({n \choose n/2+k}(n/2-k)\right)^{-1}$.

Fix $\delta>0$. Using conditions 1 and 3, we can choose $C$ such that with
probability larger than $1-\delta$, $\pi$ stays between levels
$S_{-C\sqrt{n}}$ and $S_{C\sqrt{n}}$.
It is elementary to check that there exists $C'$ (depending on $C$)
such that for all large $n$, if $k\in [-C\sqrt{n},C\sqrt{n}]$, then
$2\left({n \choose n/2+k}(n/2-k)\right)^{-1}\le \frac{C'}{\sqrt{n}2^{n}}$.

It follows that
\begin{equation} \label{eqn28}
\P(E_n\cap \pi \neq \emptyset)\le \delta+ \frac{C'|E_n|}{\sqrt{n}2^{n}},
\end{equation}
since the event $\{E_n\cap \pi \neq \emptyset\}$ is a subset of the union of the events
that $\pi$ does not stay between levels $S_{-C\sqrt{n}}$ and $S_{C\sqrt{n}}$
and that $\pi$ hits one of at most $|E_n|$ specified edges sitting between
levels $S_{-C\sqrt{n}}$ and $S_{C\sqrt{n}}$. The standard relationship between
the total influence and the edge boundary gives
$|E_n|=2^{n-1}\sum_{i\in [n]}I_i(f_n)$ and hence the last term in (\ref{eqn28})
is at most $\frac{C'\sum_{i\in [n]}I_i(f_n)}{\sqrt{n}}$. The Cauchy-Schwarz
inequality tells us this is at most $C' \sqrt{I\!I(f_n)}$ which by
assumption approaches 0 as $n\to\infty$.
Since $\delta$ is arbitrary, we are done.
\fbox{}\\

\medskip\noindent
Throughout the rest of this section we drop the superscript $n$ in the notation $\eta^n_t$.
Before starting the proof of Theorem \ref{thmmonNSimpliesCGXS},
for each $\epsilon, t>0$, we will
construct a triple coupling $(\omega,\omega^\epsilon,\eta_t)$ so that
the first two marginals have joint distribution $(\omega,\omega^\epsilon)$
(as defined right before Definition \ref{defnoise1}), and
the first and third marginals have joint distribution $(\eta_0,\eta_t),$
and with the goal that $\omega^\epsilon$ and $\eta_t$ are ``close''.
First, considering a realization of $(\omega,\omega^\epsilon),$ we let
$N^{\epsilon}_{10}=|\{x\in [n]: \omega(x)-\omega^\epsilon(x)=1\}|$
and $N_{01}^{\epsilon}=|\{x\in [n]: \omega^\epsilon(x)-\omega(x)=1\}|$.
Given $\omega,$ we have that $N^{\epsilon}_{10} \sim {\rm Bin}(|\omega|,\epsilon/2)$, $N^{\epsilon}_{01} \sim {\rm Bin}(n-|\omega|,\epsilon/2)$
and that these random variables are conditionally independent.
Similarly, given $(\eta_0,\eta_t)$
define $N_{01}^t=|\{x\in [n]: \eta_t(x)-\eta_0(x)=1\}|$. We could define
$N_{10}^t$ analogously but this random variable equals $N_{01}^t$ due to
particle conservation. The distribution of $N_{01}^t$ given $\eta_0$
is more complicated, but we will not need to know the exact distribution.

We now construct the coupling. Let $\omega=\eta_0$ be chosen uniformly.
Choose $N^{\epsilon}_{01}$, $N^{\epsilon}_{10}$ and $N^{t}_{01}$ to be conditionally
independent with their correct conditional distributions.
Given $\omega$, $N^{\epsilon}_{01}$, $N^{\epsilon}_{10}$ and $N^{t}_{01}$, we
construct a maximal coupling between $\omega^{\epsilon}$ and $\eta_t$
as follows. Assume first that
$N^{t}_{01} \leq \min(N^{\epsilon}_{01}, N^{\epsilon}_{10})$ and pick $N^{t}_{01}$ 1's of
$\omega$ uniformly and $N^{t}_{01}$ 0's of $\omega$ uniformly. Flip all of the chosen
sites and let this configuration be $\eta_t$. Continue from $\eta_t$ by picking $N^{\epsilon}_{01}- N^{t}_{01}$
($N^{\epsilon}_{10}- N^{t}_{01}$) 0's (1's) of $\omega$ uniformly among the 0's (1's)
not already picked and flip these to 1's (0's) and call the resulting configuration $\omega^{\epsilon}$.
Treat the cases $N^{\epsilon}_{10} \leq N^{t}_{01} \leq N^{\epsilon}_{01}$,
$N^{\epsilon}_{01} \leq N^{t}_{01} \leq N^{\epsilon}_{10}$ and $N^{t}_{01} \geq \max(N^{\epsilon}_{01},N^{\epsilon}_{10})$
in the obvious analogous way. This then creates a coupling
$(\omega,\omega^{\epsilon}, \eta_t)$ with the two desired 2-dimensional
marginal distributions. A {\sl crucial} feature of this coupling is that
$$
d(\omega^{\epsilon},\eta_t)=|N_{01}^{\epsilon}-N_{01}^{t}|+|N_{10}^{\epsilon}-N_{01}^{t}|.
$$
We will prove  Theorem \ref{thmmonNSimpliesCGXS}
by showing that for the correct choice of $\epsilon$,
$(\omega^{\epsilon}, \eta_t)$ satisfies conditions 1, 2 and 3 of Lemma \ref{lemma12} and then we will
apply that lemma.

\medskip

It is natural to believe that for ``most $\omega$'' $N^{\epsilon}_{10}$ and $N^{\epsilon}_{01}$ are of order
$\epsilon n/4+ O(\sqrt{n})$. We make this precise in the following lemma.
\begin{lemma} \label{lemma1}
For every $\delta>0$ and every $0<C_1<\infty$, there exists $C_2$ such that for
every $\epsilon>0$, $n$ and $\omega$ such that
$|\omega| \in [\frac{n}{2}-C_1 \sqrt{n},\frac{n}{2}+C_1 \sqrt{n}]$,
\begin{equation} \label{fluctationbound1}
\P\left(\left| N_{01}^{\epsilon}-\frac{n \epsilon}{4}\right| \geq C_2
\sqrt{n\epsilon} \ |\,  \omega\right) <\delta.
\end{equation}
\end{lemma}
\noindent
{\bf Proof.}
Fix $\delta>0$ and $C_1$.
Letting $|\omega| \in [\frac{n}{2}-C_1 \sqrt{n},\frac{n}{2}+C_1 \sqrt{n}]$, obviously
\begin{eqnarray*}
\lefteqn{ \left| N_{01}^{\epsilon}-\frac{n \epsilon}{4}\right|}\\
& &  \leq
\left| N_{01}^{\epsilon}-\E[N_{01}^{\epsilon} |\,  \omega] \right|
+\left| \E[N_{01}^{\epsilon} |\,  \omega]-\frac{n \epsilon}{4}\right| \\
& & \leq
\left| N_{01}^{\epsilon}-\E[N_{01}^{\epsilon} |\,  \omega] \right|
+\frac{\epsilon}{2}C_1 \sqrt{n}.
\end{eqnarray*}
Therefore it suffices to show that there exists a $C$ such that
for every $\epsilon>0$ and $n$,
\[
\P\left(\left| N_{01}^{\epsilon}-\E[N_{01}^{\epsilon}|\,  \omega]\right| \geq
C \sqrt{n\epsilon} \ |\,  \omega\right) <\delta.
\]
This follows in fact for any $\omega$ by Markov's inequality since
\begin{eqnarray*}
\lefteqn{\P\left(\left| N_{01}^{\epsilon}-\E[N_{01}^{\epsilon} |\,  \omega]\right| \geq C \sqrt{n\epsilon} \ |\,  \omega\right) }\\
& & \leq \frac{{\rm Var}(N_{01}^{\epsilon} |\,  \omega)}{C^2 n\epsilon}
=\frac{\epsilon/2 (1-\epsilon/2) (n-|\omega|)}{C^2 n\epsilon}
\leq \frac{1}{C^2}<\delta,
\end{eqnarray*}
by taking $C$ large enough.
\fbox{}\\

\medskip


We now proceed by proving the corresponding result for $N_{01}^t$.
It will be useful to write
\[
N_{01}^t=\sum_{x:\eta_0(x)=0} I_{\eta_t(x)=1}.
\]
Given $\eta_0$ and $x$ with $\eta_0(x)=0$, the event that
$\eta_t(x)=1$ is the event that the Poisson clock corresponding to some edge
with $x$ as an endpoint rings before time $t$ and the last time such a
clock rings $x$ was updated to a 1.
Hence, if $\eta_0(x)=0$, then
$\P(\eta_t(x)=1 |\,  \eta_0)=(1-e^{-t})|\eta_0|/n$
and so
\begin{equation} \label{eqn24}
\E[N_{01}^t |\,  \eta_0]=(1-e^{-t}) \frac{|\eta_0|(n-|\eta_0|)}{n}.
\end{equation}
Furthermore we have the following lemma.

\begin{lemma} \label{lemma2}
For every $\delta>0$ and every $0<C_1<\infty$, there exists $C_3$ such that for
every $t>0$, $n$ and $\eta_0$ such that
$|\eta_0| \in [\frac{n}{2}-C_1 \sqrt{n},\frac{n}{2}+C_1 \sqrt{n}]$,
\begin{equation} \label{fluctationbound2}
\P\left(\left| N_{01}^{t}-\frac{(1-e^{-t})n}{4} \right| \geq C_3 \sqrt{n}
\sqrt{1-e^{-t}}\ |\,  \eta_0\right) <\delta.
\end{equation}
\end{lemma}
\noindent
{\bf Proof.}
Fix $\delta>0$ and $0<C_1<\infty$. Then for all $t>0$, $n$ and $\eta_0$ satisfying the
conditions in the statement, by the triangle inequality, (\ref{eqn24}) and the assumption
on $\eta_0$, we have
\begin{eqnarray*}
\lefteqn{\left| N_{01}^{t}-\frac{(1-e^{-t})n}{4} \right|}\\
& & \leq \left| N_{01}^{t}-\E[N_{01}^t |\,  \eta_0] \right|
+\left| \E[N_{01}^t |\,  \eta_0]-\frac{(1-e^{-t})n}{4} \right| \\
& & \leq \left| N_{01}^{t}-\E[N_{01}^t |\,  \eta_0] \right|+C_1^2(1-e^{-t}).
\end{eqnarray*}

Therefore it suffices to show that there exists a $C$ (depending on $\delta$ and $C_1$)
such that for every $t>0$, $n$ and $\eta_0$ satisfying the conditions in the statement,
we have that
\begin{equation} \label{eqn25}
\P\left(\left| N_{01}^{t}-\E[N_{01}^t |\,  \eta_0] \right| \geq C \sqrt{n}
\sqrt{1-e^{-t}} \ |\,  \eta_0\right) <\delta.
\end{equation}
This in fact holds for all $\eta_0$ as we shall see. By Markov's inequality
\begin{equation} \label{eqn26}
\P\left(\left| N_{01}^{t}-\E[N_{01}^t |\,  \eta_0] \right| \geq C \sqrt{n}
\sqrt{1-e^{-t}}\ |\,  \eta_0\right)
\leq \frac{{\rm Var}(N_{01}^{t}|\,  \eta_0 )}{C^2 n(1-e^{-t})},
\end{equation}
and therefore we are interested in bounding ${\rm Var}(N_{01}^{t}|\,  \eta_0).$
Consider therefore $\E[(N_{01}^{t})^2 |\,  \eta_0 ]$. Given $\eta_0$ and $x\neq y$
such that $\eta_0(x)=\eta_0(y)=0$, it is easy to see that the events
$\{\eta_t(x)=1\}$ and $\{\eta_t(y)=1\}$ are negatively correlated
and therefore for such an $\eta_0$
\begin{eqnarray*}
\lefteqn{\P(\eta_t(x)=\eta_t(y)=1 |\,  \eta_0)}\\
& & \leq \P(\eta_t(x)=1 |\,  \eta_0)^2
=\left(\frac{(1-e^{-t})|\eta_0|}{n}\right)^2.
\end{eqnarray*}
Therefore
\begin{eqnarray*}
\lefteqn{\E[(N_{01}^{t})^2 |\,  \eta_0 ]
= \E\left[\left(\sum_{x:\eta_0(x)=0} I_{\eta_t(x)=1}\right)^2 |\,  \eta_0 \right]}\\
& & \leq (1-e^{-t}) \frac{|\eta_0|(n-|\eta_0|)}{n}
+(n-|\eta_0|)(n-|\eta_0|-1)\left(\frac{(1-e^{-t})|\eta_0|}{n}\right)^2.
\end{eqnarray*}
Hence by (\ref{eqn24}),
\[
{\rm Var}(N_{01}^{t}|\,  \eta_0 ) \leq n(1-e^{-t}),
\]
and so (\ref{eqn25}) follows by choosing $C$ large enough in (\ref{eqn26}).
\fbox{}\\

\medskip\noindent
For the proof of Theorem \ref{thmmonNSimpliesCGXS}, we will
only need weaker versions of Lemmas \ref{lemma1} and \ref{lemma2}
where the $\sqrt{\epsilon}$ term in (\ref{fluctationbound1}) and the
$\sqrt{1-e^{-t}}$ term in (\ref{fluctationbound2}) do not appear. However,
for the proof of Theorem \ref{thmmonCGXStableimpliesNstable}, the stronger
versions which are stated will be needed.

We can now easily prove Theorem \ref{thmmonNSimpliesCGXS}.

\noindent
{\bf Proof of Theorem \ref{thmmonNSimpliesCGXS}.}
Fix $t>0$ and let $\epsilon=1-e^{-t}$. Consider the
triple $(\omega, \omega^{\epsilon},\eta_t)$ with
$\omega=\eta_0$ defined earlier.
By Theorem \ref{BKSMAIN} \cite{BKS}, the assumption that
$\lim_{n \rightarrow \infty}I\!I(f_n)=0$ implies that
$\{f_n\}_{n \in \N}$ is NS and since
\begin{eqnarray} \label{eqn:starstar}
\lefteqn{| \E[f_n(\omega)f_n(\omega^{\epsilon})]-\E[f_n(\eta_0)f_n(\eta_t)] |} \\
& & \label{eqn:starstar1}
=2|\P(f_n(\omega)=f_n(\omega^{\epsilon}))-\P(f_n(\eta_0)=f_n(\eta_t))|  \nonumber\\
& & \label{eqn:starstar2} \leq 2\P(f_n(\omega^{\epsilon}) \neq f_n(\eta_t)), \nonumber
\end{eqnarray}
it suffices to show that
\[
\lim_{n \rightarrow  \infty} \P(f_n(\omega^{\epsilon})\neq f_n(\eta_t))=0.
\]

We do this by showing that $(\omega^{\epsilon},\eta_t)$ satisfies conditions 1, 2 and 3
of Lemma \ref{lemma12}. Conditions 1 and 2 are immediate from the construction of
$(\omega^{\epsilon},\eta_t)$. Furthermore, by our choice of $\epsilon$
and our coupling, we have
\begin{eqnarray*}
\lefteqn{d(\omega^{\epsilon},\eta_t)=|N_{01}^{\epsilon}-N_{01}^{t}|+|N_{10}^{\epsilon}-N_{01}^{t}|}\\
& & \leq \left|N_{01}^{\epsilon}-\frac{\epsilon n}{4}\right|+\left|N_{10}^{\epsilon}-\frac{\epsilon n}{4}\right|
+2 \left|N_{01}^{t}-\frac{(1-e^{-t}) n}{4}\right|.
\end{eqnarray*}
For fixed $\delta>0$, choose $C_1$ such that
\[
\P(|\omega| \in [n/2-C_1 \sqrt{n}, n/2+C_1 \sqrt{n}]) >1-\delta/2.
\]
Then use Lemmas \ref{lemma1} and \ref{lemma2} to find $C$ such that for
any $n$ and for any $\omega$ with
$|\omega| \in [n/2-C_1 \sqrt{n}, n/2+C_1 \sqrt{n}]$,
\begin{eqnarray*}
\P\left(\left|N_{01}^{\epsilon}-\frac{\epsilon n}{4}\right|+\left|N_{10}^{\epsilon}-\frac{\epsilon n}{4}\right|
+2 \left|N_{01}^{t}-\frac{(1-e^{-t}) n}{4}\right| \leq C \sqrt{n} |\,  \omega \right) \\
>1-\delta/2.
\end{eqnarray*}
(This would still be true if $\sqrt{n}$ were replaced by $\sqrt{n\epsilon}$.)
This verifies condition 3 and so Lemma \ref{lemma12} proves the theorem.
\fbox{}\\

\noindent
{\bf Remark:} We quickly give a sketch of a different proof of
Corollary \ref{cor:monotone} which does not go through the use of
Lemma \ref{lemma12}. To do this, first it is not so difficult to go
from the 3-way
coupling already described and prove the following claim.
For all $\epsilon >0$ and $\delta >0$, there exists $k$ such that for all
$n$ large, there exists a 5-way coupling
$(\omega,\omega^\epsilon,\eta_t,Y,W)$ such that:
\begin{enumerate}
\item $(\omega,\omega^\epsilon,\eta_t)$ has the distribution given earlier.
\item $(\omega^\epsilon,W,Y)$ has distribution
$\pi_{1/2}\times\pi_{1-k/\sqrt{n}}\times\pi_{k/\sqrt{n}}$
\item $\P(\omega^\epsilon\wedge W \le \eta_t\le \omega^\epsilon\vee Y) \ge 1-\delta$.
\end{enumerate}
Next, one can show that for noise sensitive sequences, the ``threshold'' window
is of size much larger than $1/\sqrt{n}$; this is due to the fact
that they are uncorrelated with majority-type functions. It follows that
$f(\omega^\epsilon\wedge W) =f(\omega^\epsilon\vee Y)$ with probability
going to 1. Since on most of the space,
$\eta_t$ is trapped between these two configurations, it follows by monotonicity
that $f(\omega^\epsilon) =f(\eta_t)$ with probability going to 1.

\bigskip\noindent
We now show that Proposition \ref{prop:paradox} follows from
Lemma \ref{lemma12}.

\medskip\noindent
{\bf Proof of Proposition \ref{prop:paradox}.}
Apply Lemma \ref{lemma12} with $(\xi_n,\xi^*_n)$ equal to
$(\omega,\omega^{\epsilon_n})$. It is easy to check
that all the conditions of the lemma are met.
\fbox{}\\

\noindent Corollary \ref{corr1} follows easily from the proof of
Theorem \ref{thmmonNSimpliesCGXS}.

\medskip\noindent{\bf Proof of Corollary \ref{corr1}.}
Using the above coupling as in Theorem \ref{thmmonNSimpliesCGXS}, we see as in the proof of
Theorem \ref{thmmonNSimpliesCGXS} that the pair $(\omega^{\eps_n},\eta_{t_n})$
satisfies conditions 1,2 and 3 of Lemma \ref{lemma12}.
This lemma therefore implies
that $\P(f_n(\omega^{\epsilon_n}) \neq f_n(\eta_{t_n}))\to 0$ and hence
(\ref{eqn:starstar})-(\ref{eqn:starstar2}) gives the result.
\fbox{}\\

\medskip\noindent
We now turn to the proof of Theorem \ref{thmmonCGXStableimpliesNstable},
where we will use the ideas from the proof of Theorem \ref{thmmonNSimpliesCGXS}.
We first have the following variant of Lemma \ref{lemma12}.

\begin{lemma}\label{lemma12newversion}
There exists a universal constant $\Lambda$ so that for all $\delta>0$,
we have the following for large $n$. If
$(\omega,\omega^*)$ is a pair of random
configurations in $\{0,1\}^{[n]}$ satisfying the following conditions:

\begin{enumerate}
\item $\omega$ is uniformly distributed on  $\{0,1\}^{[n]}$

\item $\omega \preceq \omega^*$

\item The joint distribution of $(\omega,\omega^*)$ is invariant under all
permutations of $[n]$

\item
\[
d(\omega,\omega^*)\le \delta \sqrt{n} \,\, a.s.,
\]

\end{enumerate}
then, for any monotone Boolean function $f$ on $\{0,1\}^{[n]}$,
\[
\P(f(\omega)\neq f(\omega^*))\le \Lambda\delta.
\]
\end{lemma}

\noindent
{\bf Proof.}
We may assume that $f$ maps into $\{0,1\}$. Fix $\delta>0$.
Let $\sigma\in\{0,1\}^{[n]}$ have distribution $\pi_{1/2+2\delta/\sqrt{n}}$. We claim that
there exists a coupling $(\omega,\omega^*,\sigma)$ so that
$(\omega,\omega^*)$ has the correct conditional distribution,
$\omega\preceq \sigma$ a.s.\ and
$\omega^*\preceq \sigma$ with probability going to 1 as $n\to\infty$.
To see this, one first couples
$\omega$ and $\sigma$ in the canonical (monotone) way.
One notes that the number of coordinates which are 0 for $\omega$ and
1 for $\sigma$ has a Binomial distribution with parameters $n$ and
$2\delta/\sqrt{n}$. Chebyshev's inequality tells us that such a random variable
is larger than $\delta\sqrt{n}$ with probability going to 1 as $n\to\infty$.
If $|\sigma|-|\omega|\ge \delta\sqrt{n}$, we let $\omega^*$ be
a configuration between $\omega$ and $\sigma$ such that $|\omega^*|$ is chosen
with the correct conditional (given $\omega$) distribution and then chosen uniformly
among those between $\omega$ and $\sigma$.
If $|\sigma|-|\omega|< \delta\sqrt{n}$, just choose
$\omega^*$  with the correct conditional (given $\omega$) distribution.
This proves the claimed coupling.

By monotonicity, it now suffices to show that
$\P(f(\omega)\neq f(\sigma))\le \Lambda\delta$. It is clear that
this probability is
$$
\E_{\pi_{1/2+2\delta/\sqrt{n}}}[f]-\E_{\pi_{1/2}}[f].
$$
Letting $g(p)$ denote the expected value of $f$ under $\pi_p$, it
suffices to show that the derivative of $g(p)$ in the interval $[1/2,3/4]$
is bounded by $\sqrt{n}$ times a universal constant for any monotone
function. This however is easy and well known.
\fbox{}

\medskip\noindent
{\bf Proof of Theorem \ref{thmmonCGXStableimpliesNstable}.}
Let $\epsilon$ and $t$ be related by $\epsilon = 1-e^{-t}$.
We consider again the triple coupling $(\omega, \omega^{\epsilon},\eta_t)$ that
we defined earlier in this section.
As in the proof of Theorem \ref{thmmonNSimpliesCGXS},
Lemmas \ref{lemma1} and \ref{lemma2}
yield that for any $\alpha>0$, there exists $C=C(\alpha)$
such that for all $\epsilon$ and $n$,
$$
\P(d(\omega^{\epsilon},\eta_t)> C \sqrt{\epsilon} \sqrt{n}) < \alpha.
$$
To prove NStable, fix $\delta>0$ arbitrarily. Then let $\alpha=\delta/10$ and
$C=C(\alpha)=C(\delta/10)$. Since our sequence is
CGXStable, we know that for $t$ sufficiently small,
$\P(f_n(\omega) \neq f_n(\eta_t))<\delta/10$ for all $n$. Furthermore,
choose t sufficiently small so that
\begin{equation} \label{eqn3}
\epsilon= 1-e^{-t} \le \frac{\delta^2}{C^2}.
\end{equation}
Now for this choice of $\epsilon$, we have that for every $n$,
$$
\P(f_n(\omega)\neq f_n(\omega^\epsilon)) \le
\P(f_n(\omega)\neq f_n(\eta_t)) +
\P(d(\omega^{\epsilon},\eta_t)>C \sqrt{\epsilon} \sqrt{n})
$$
$$
+
\P(f_n(\omega^\epsilon)\neq f_n(\eta_t),
d(\omega^{\epsilon},\eta_t)\le C \sqrt{\epsilon} \sqrt{n}).
$$
Each of the first two terms are at most $\delta/10$ while (\ref{eqn3})
tells us that the last term is at most
$$
\P(f_n(\omega^\epsilon)\neq f_n(\eta_t),
d(\omega^{\epsilon},\eta_t)\le \delta \sqrt{n}).
$$
Letting $\xi$ be $\max\{\omega^\epsilon,\eta_t\}$,
the last term is then at most
$$
\P(f_n(\omega^\epsilon)\neq f_n(\xi),
d(\omega^{\epsilon},\xi)\le \delta \sqrt{n})+
\P(f_n(\eta_t)\neq f_n(\xi),
d(\eta_t,\xi)\le \delta \sqrt{n}).
$$
We bound only the first summand, the second is handled identically.
Defining $\tilde{\xi}$ to be $\xi$ on
$d(\omega^{\epsilon},\xi)\le \delta\sqrt{n}$ and
$\omega^{\epsilon}$ otherwise, the first term becomes
$$
\P(f_n(\omega^\epsilon)\neq f_n(\tilde{\xi})).
$$
Lemma \ref{lemma12newversion} tells us that for fixed $\delta$, this is
at most $\Lambda\delta$ for large $n$. Since
$\Lambda$ is a fixed constant and we can take $\epsilon$ even smaller
to insure what we want for finitely many values of $n$,
NStability is established. \fbox{}\\

\section{Spectral techniques and percolation} \label{secthm5}

As pointed out in the introduction,
Theorem \ref{thm5} follows from Theorem \ref{thmPCNS} and
Corollary \ref{cor:monotone}.  However, we give an alternative
proof which uses the geometric description of the spectral sample,
does not use monotonicity and will be applicable to the medium
range exclusion process.
We also believe this approach is potentially applicable in other situations,
for instance when the functions are non-monotone.
We start with a lemma.

\begin{lemma} \label{lemmaAB}
Let $Q_n$ be a sequence of finite sets, and for every $n$ let $\nu_n$ be a
subprobability measure on $Q_n$. Furthermore, let $P_n$ be a symmetric
transition matrix on the set $Q_n$ and assume that there exists a sequence
$B_n$ of subsets of $Q_n$ such that
$\nu_n(B_n^c)\to 0$ and $\max_{x\in \Lambda^2_n}P_n(x,\Lambda^2_n)\to 0$.
Then
$$
\sum_{x,y\in Q_n}  \sqrt{\nu_n(x)}\sqrt{\nu_n(y)}P_n(x,y)\to 0.
$$
\end{lemma}

\noindent
{\bf Proof:} Fix  $\epsilon >0$. Let $\Omega_n=\{(x,y)\in Q_n\times Q_n:
\sqrt{\nu_n(x)}\sqrt{\nu_n(y)}  \le \epsilon(\nu_n(x)+ \nu_n(y))\}$.
We break the sum into four sums. The first is over $\Omega_n$,
the second is over $\Omega_n^c\cap \{y\not\in B_n\}$,
the third is over $\Omega_n^c\cap \{y\in B_n\}\cap \{x\not\in B_n\}$
and the fourth is over $\Omega_n^c\cap \{y\in B_n\}\cap \{x\in B_n\}$.

For the first sum, the summand can be replaced by
$\epsilon(\nu_n(x)+ \nu_n(y))P_n(x,y)$. Summing this over all $x$ and $y$
and using the fact that $P_n$ is a symmetric transition matrix gives
an upper bound of $2\epsilon$.

Next, it is easy to check that if $(x,y)$ is not in $\Omega_n$, then the ratio
of $\max\{\sqrt{\nu_n(x)},\sqrt{\nu_n(y)}\}$
and $\min\{\sqrt{\nu_n(x)},\sqrt{\nu_n(y)}\}$ is at most $1/\epsilon$.
This will be used in the last three sums.

For the second sum, the summand can be bounded above by $\nu_n(y)
P_n(x,y)1/\epsilon$. When we sum first over $x$ and then over
$y$, we get $\nu_n(B_n^c)/\epsilon$.

For the third and fourth sums, the summand can be upper bounded by $\nu_n(x)
P_n(x,y)1/\epsilon$. The third sum then becomes at most
$\nu_n(B_n^c)/\epsilon$ while the fourth sum becomes at most
$\max_{x\in B_n}P(x,B_n)/\epsilon$.

Our assumptions then give that the limsup of the full double
sum is at most $2\epsilon$. Since $\epsilon$ is arbitrary, we are done.
\fbox{}\\

We have the following theorem which easily follows from Lemma \ref{lemmaAB}.
Let $\calP(U)$ denote the set of subsets of $U$, and for $A_n \in \calP(V_n),$
we extend notation and write
\[
P^n_t(S, A_n):=\sum_{S'\in A_n}P^n_t(S, S').
\]

\begin{theorem} \label{thm3}  Consider the exclusion process
with respect to $\{G_n, \alpha_n\}_{n\in\N}$.
Let ${\cal S}_n$ denote the spectral sample of $f_n$.
If there exist $t>0$ and sets $A_n\subseteq \calP(V_n)\setminus \{\emptyset\}$ such that
$\P({\cal S}_n\in (A_n\cup\{\emptyset\})^c)\rightarrow 0$ and
$\max_{S\in A_n} P^n_t(S, A_n)\rightarrow 0$, then $\{f_n\}_{n\in \N}$ is XS.
\end{theorem}
\noindent
{\bf Proof.}
In Lemma \ref{lemmaAB}, let $Q_n$ be $\calP(V_n)$, $\nu_n$ be the spectral
measure of $f_n$ restricted to $Q_n\setminus \{\emptyset\}$
(i.e. $\nu_n(S)=\hat{f}_n(S)^2$ for every $S\neq \emptyset$ and
$\nu_n(\emptyset)=0$), $P_n$ be $P^n_t$ and $B_n$ be $A_n$.
Lemma \ref{lemmaAB} and Proposition \ref{prop2} now imply XS.
\fbox{}\\

\noindent
{\bf Remarks:}\\
1. In fact Lemma \ref{lemmaAB} implies that even (\ref{eqn4}) holds.\\
2. Observe that the main assumption of Theorem \ref{thm3} is essentially that
${\cal S}_n$ and ${\cal S}_nP^n_t$ are very singular off of
${\cal S}_n=\emptyset$.

\medskip\noindent
As an application of Theorem \ref{thm3} we will now prove Theorem \ref{thm5}.
We point out that from this point on in the paper,
we will rely on some nontrivial results from \cite{GPS}. Specific pointers to \cite{GPS} will be given when needed.

\bigskip
\noindent
{\bf Proof of Theorem \ref{thm5}.} \\
We prove the statement for $a=b=1,$ and for simplicity we shift the boxes so that
$B_n:=[-n/2,n/2]^2\cap \T$ becomes the relevant box. In addition, we will indicate
the changes needed for the general case $a,b>0.$ For $0<\beta<1$, let
$B_{n^{1-\beta}}:=[-n^{1-\beta}/2,n^{1-\beta}/2]^2\cap \T$.
Let $f_n=2I_{C_{n}}-1$ where $C_{n}$ is as in the statement of
Theorem \ref{thm5} (although now defined on $B_n$). Let ${\cal S}_n$ be the spectral sample of $f_n$. It follows from
Lemma 3.2 in \cite{GPS} that
\begin{equation} \label{eqn7}
\P({\mathcal S}_n \cap B_{n^{1-\beta}}\neq \emptyset)\leq n^{-(5/4+o(1))\beta}.
\end{equation}
Note next that obviously
\[
\P({\mathcal S}_n \cap B_{n^{1-\beta}}\neq \emptyset)=\sum_{S:S\cap B_{n^{1-\beta}}\neq \emptyset}\hat{f_n}(S)^2.
\]
Furthermore, due to 
Theorem 7.4 in \cite{GPS}, we have that for any $\epsilon>0$,
\[
\lim_{n\to\infty}\P(n^{3/4-\epsilon}\leq |{\mathcal S}_n|\leq n^{3/4+\epsilon})= 1.
\]
(If we used a rectangle which is not a square, then $\E[f_n]\not\to 0$ and we
would need to include the event $\{{\mathcal S}_n=\emptyset\}$ in the above to
get a probability going to 1.) We will apply Theorem \ref{thm3} to the sets
\[
A_n:=\{S\subseteq B_n:S\cap B_{n^{1-\beta}}= \emptyset, n^{.74}\leq |S|\leq n^{.76}\}.
\]
By the above $\P(S_n\in (A_n\cup \emptyset)^c)\rightarrow 0$ as $n\to\infty$.
If we can show that
$\max_{S\in A_n} P^n_t(S,A_n)\rightarrow 0$ for any $t>0$, then Theorem \ref{thm3}
implies the result.

Fix $t>0$. For every $n$, with $S\in A_n$, write $S^n_t$ for the random set $S_t$
that appears in Lemma \ref{lemma3}.
Using Lemma \ref{lemma3} we conclude that for any $\delta>0$,
$\P(|S^n_t\cap S^c|\geq n^{.73} )\geq 1-\delta$ for all $n$ large enough,
uniformly in $S\in A_n$.
Since the set $S^n_t\cap S^c$ is distributed uniformly at random on
$S^c$ we have that for any $\beta$
\[
\P(S^n_t\cap S^c \cap B_{n^{1-\beta}}=\emptyset |\, |S^n_t\cap S^c|\geq n^{.73})\leq \left(1-n^{-2\beta} \right)^{n^{.73}}
\]
and so
\[
\P(S^n_t\cap S^c \cap B_{n^{1-\beta}}=\emptyset)\leq \delta +\left(1-n^{-2\beta} \right)^{n^{.73}}
\]
for all $S\in A_n$ if $n$ is large.
Therefore, if $\beta\le .36$ we can let $n\to\infty$ first and then
$\delta\to 0$ to conclude that
\[
\lim_{n\to\infty}\P(S^n_t\cap S^c \cap B_{n^{1-\beta}}=\emptyset) = 0,
\]
uniformly in $S\in A_n$. We conclude that
\[
\max_{S\in A_n}P^n_t(S,A_n) \leq \max_{S\in A_n}\P(S^n_t\cap S^c \cap B_{n^{1-\beta}}=\emptyset)\rightarrow 0.
\]
\fbox{}\\

\section{Medium-range dynamics for percolation}\label{s.shortrange}
We will start this section by proving Proposition \ref{pr.SR2}, a weaker version of Theorem \ref{th.SR}.
This has the advantage of conveying the main ideas while at the same time being useful in the proof
of Theorem \ref{th.SR}.
\begin{proposition}\label{pr.SR2}
For any $\alpha>3/8$ and $t>0$, letting $f_n=I_{C_n}$, we have
\[
\lim_{n \to \infty}\E[f_n(\eta^{\alpha,n}_0)f_n(\eta^{\alpha,n}_t)]-\E[f_n(\eta^{\alpha,n}_0)]^2 = 0.
\]
\end{proposition}

\noindent
Intuitively, the idea behind the proof can be explained as follows.
We will consider a partition of $[0, an] \times [0, b n]$ into squares
of sidelength $r<<n^{\alpha}$. We proceed to show
that with high probability no such $r-$square receives ``too many points'' from ${\cal S}$ after running the dynamics,
and that a certain fraction of the points in ${\cal S}$ must move.
Combining these two results we will conclude that there must be ``many''
$r-$squares that receive at least one point from ${\cal S}$.
The new configuration will then be ``singular to ${\cal S}$''
and we will be able to apply Theorem \ref{thm3}.

\medskip

Before we can present the proof of Proposition \ref{pr.SR2}, we will need some
intermediate results. First we need a definition.

\begin{definition} \label{defn:posdef}
A symmetric function $f : V^k \to \R$ is called
\emph{positive definite} if for any choice of $x_1, \ldots, x_k \in V$,
and any $1\le i<j\le k$, the ``matrix''
\[
k^{i,j}(x,y):= f(x_1, .., x_{i-1}, x , .. , y,x_{j+1}, ...,x_k)
\]
is positive definite.
\end{definition}

\begin{proposition}\label{propIPS}[\cite{L} Proposition 1.7, chapter 8]
Let $G=(V,E)$ be a graph and fix any initial set
$A:= \{z_1, \ldots, z_k\}\subseteq V$.
Let $\hat X_t=\{\hat X^1_t,\ldots, \hat X^k_t\}$ denote any symmetric
exclusion process on $G$ starting from $A$ and let
$X_t=\{X^1_t, \ldots, X^k_t\}$
denote the trajectory of $k$ independent random walks starting from $A$,
each evolving independently in the same way as one particle in the
exclusion process.
If $f: V^k \to \R$ is a symmetric positive definite function, then
\[
\Eb{f(\hat X_t)} \le \Eb{f(X_t)}.
\]
\end{proposition}

\medskip \noindent
This will allow us to prove the following lemma where $\P^{\alpha,n}$
refers to the process $\{\eta^{\alpha,n}_t\}$ and
$A \mapsto_t B$ denotes the event that the set $A$ has been moved to the set
$B$ under the corresponding random permutation $\pi_t$.

\begin{lemma} \label{lemma:bound.on.points.going.to.points}
For any $\alpha>0$, n, $t>0$, and $2k$ distinct points
$z_1,\ldots, z_k,x_1,\ldots, x_k$,
\[
\P^{\alpha,n}[ z_1 \mapsto_t x_1, \ldots , z_k \mapsto_t x_k] \le  \Bigl( \frac k {n^{2\alpha}} \Bigr)^k\,.
\]
\end{lemma}

\noindent
{\bf Proof.}
Apply Proposition \ref{propIPS} to
\[
f(y_1,\ldots, y_k):= 1_{\{ y_1,\ldots, y_k\} \subseteq \{x_1,\ldots, x_k\} }\,.
\]
The function $f$ clearly satisfies the conditions and therefore one obtains that
\begin{eqnarray*}
\lefteqn{\P^{\alpha,n}[ z_1 \mapsto_t x_1, \ldots , z_k \mapsto_t x_k] \leq
\P^{\alpha,n}( \{z_1,\ldots,z_k\} \mapsto_t \{x_1, \ldots, x_k\})}\\
& & = \E^{A=\{z_1,\ldots, z_k\}}  [f(\hat X_t)] \le \E^{A=\{z_1,\ldots, z_k\}} [f(X_t)] \\
& &= \P^{A=\{z_1,\ldots, z_k\}} [X_t \subseteq \{ x_1, \ldots, x_k \}] \le  \Bigl( \frac k {n^{2\alpha}} \Bigr)^k \,.
\end{eqnarray*}
The last inequality follows since a fixed particle $z_j$ must in its last jump
land in a set of size $k$ while the range of a step is $n^{2\alpha}$ ; note that the
assumption of the disjointness of our $2k$ points is used here to insure that there is 
some jump associated to each of $z_1,\ldots,z_k$.
This technique (a kind of symmetrization) is used for example in [\cite{JL},
Section 6].
\fbox{}

We can now prove the following proposition.
\begin{proposition}\label{pr.NB}
Given two sets $E$ and $F$, let $N_{t}(E,F)$ be the number of
points in $E$ at $t=0$ which are in $F$ at time $t$.
Then, for all integers $k\ge 1$, there exists an absolute constant $C_k>0$ such that
for any disjoint sets $E$, $F$, and for any $\alpha>0,n$ and $t>0$,
\begin{equation}\label{e.NB}
\E^{\alpha,n}[N_{t}(E,F)^k] \le C_k \,  \Bigl( 1 \vee \frac {|E||F|}{n^{2\alpha}} \Bigr)^k\,.
\end{equation}
\end{proposition}
\noindent
{\bf Proof.}
Obviously
\begin{eqnarray}
\E^{\alpha,n}[N_{t}(E,F)^k]
& = & \sum_{z_1, \ldots, z_k \in E} \sum_{x_1,\ldots, x_k \in F}
\P^{\alpha,n}[z_1 \mapsto_t x_1, \ldots , z_k \mapsto_t x_k,
\end{eqnarray}
where the $z_i$'s and $x_i$'s appearing need not be distinct.
We first consider the non-diagonal terms
(i.e. where all the points $z_i$ and hence the
$x_i$'s are distinct),
since the non-diagonal terms involve lower moments.
There are at most $|F|^k$ possible choices for the
non-diagonal entries $\{x_1, \ldots, x_k \}$ and at
most $|E|^k$ choices for the set of initial points
$\{z_1, \ldots, z_k \}$. Using Lemma
\ref{lemma:bound.on.points.going.to.points}, we see that
the sum of the non-diagonal terms is bounded by
$C_k \Bigl(\frac {|E||F|}{n^{2\alpha}} \Bigr)^k$, where $C_k$
depends on $k$ but is independent of $n.$ Diagonal terms,
which correspond to the non-diagonal terms for lower moments,
are treated the same way. Were it the case that
$\frac {|E||F|}{n^{2\alpha}}<1$,
then the diagonal terms would actually dominate which is why we need to take
the maximum with 1 in \eqref{e.NB}.
\fbox{}

\medskip \noindent
{\bf Proof of Proposition \ref{pr.SR2}.}
Fix $\alpha > 3/8$ and $t>0$.
Let $1\ll r_n \ll n$ be a mesoscopic scale that will grow with $n$ which
we will choose later and
partition $[0,a n]\times[0,bn]$ into squares of sidelength $r_n$. (Here we
ignore the fact that $an$ and/or $bn$ might not be divisible by $r_n$ since this
is easily handled.)
Observe that there are $O(n^2/r_n^2)$ such squares.

Fix $\eps>0$ to be determined later and let
\begin{eqnarray*}
\lefteqn{A_n=A_n^{r_n,\eps} := \{S \subseteq [0,an]\times[0,bn]:} \\
& & n^{3/4-\eps}<|S| < n^{3/4+\eps},
 \,\, |S_{(r_n)}| \le (n/r_n)^{3/4+\eps} \}\,,
\end{eqnarray*}
where $S_{(r_n)}$ is the set of $r_n$-squares in
$[0,a n]\times[0,bn]$ which intersect $S$.

Let $\Spec_n$ have the distribution of the spectral distribution of $f_n$. It follows
from Lemma 3.2 in \cite{GPS} as well as the derivation of the identity (7.3) in \cite{GPS} that the expected number of $r_n$-squares in $[0, an] \times [0, bn]$ intersected by $\Spec_n$ 
is up to constants $\frac {n^2} {r_n^2} \alpha_4(r_n,n)$. In particular, by Markov, with high probability the number of
$r_n$-squares in $[0,a n]\times[0,bn]$ intersected by $\Spec_n$ is at
most $(n/r_n)^{3/4+o(1)}$. Furthermore, it also follows that there exists a sequence
$\{\delta_n\}_{n\in \N}$ such that $\delta_n>0$ for every $n,$ $\delta_n \to 0$
and on the event
$\Spec_n\neq \emptyset$, the conditional probability that
$|\Spec_n|$ is larger than $n^{3/4-\delta_n}$ tends to 1. It follows that
for any $r_n \le n^{1-\eps}$ we have
\[
\Pb{ \Spec_n \in (A_n^{r_n,\eps}\cup \emptyset)^c }\to 0\,,
\]
as $n\to \infty$. Using a modified version of Theorem \ref{thm3}
(which is proved analogously and left to the reader),
it remains to show, for a well chosen value of $r_n$, that
for our medium-range dynamics:
\[
\max_{S\in A_n}P^{\alpha,n}_t(S,A_n) \rightarrow 0\,,
\]
as $n\to \infty$.

Assume for the moment that $r_n$ is chosen so that
$n^{3/4-\eps}\frac{r_n^2}{n^{2\alpha}} \ge 1$ for large $n$. Fix
$S\in A_n$. For $B$ an $r_n$-box, consider the random variable
$N_t(S\cap B^c,B)$ (note that it depends on $\alpha$ and $n$) defined in Proposition \ref{pr.NB}.
By our choice of $r_n$, Proposition \ref{pr.NB} tells us that the $k$th
moment of this random variable is, for large $n$, at most
$C_k\Bigl(|S| \frac {r_n^2}{n^{2\alpha}} \Bigr)^k$. This
implies that for any $r_n$-square $B$ in $[0,an]\times [0,bn]$, we have that
with high probability $B$ will receive \emph{few points}
from $S\cap B^c$. Indeed, one has for large $n$:
\[
\P^{\alpha,n}\left(N_t(S\cap B^c,B) > \left(\frac{n}{r_n}\right)^\eps |S| \frac{r_n^2}{n^{2\alpha}}\right) \le C_k \left(\frac{r_n}{n}\right)^{k \eps}.
\]

Let us choose $k=k_\eps$, such that $k \eps >2$. Since there are $O(n^2/r_n^2)$
$r_n$-squares in $[0,an]\times[0,bn]$ and since furthermore we assumed
$r_n \le n^{1-\eps}$, we conclude that if $H_n$ is the event (depending on $S$) that all
$r_n$-squares $B$ in our domain are such that
$N_t(S\cap B^c,B) \le (\frac{n}{r_n})^\eps |S| \frac{r_n^2}{n^{2\alpha}}$,
then $\min_{S\in A_n}\Pb{H_n} \to 1$ as $n \to \infty$. 

We now want a uniform (in $S$) upper bound on $P^n_t(S,A_n)$. Let $K_n$ be
the event that at least $n^{3/4-2\eps}$ points in $S$ are in a different
$r_n$-square of $[0,an]\times[0,bn]$
at time $t$ than they were at time $0$. Clearly, if we
assume that $r_n= o(n^{\alpha})$ (here $n^{\alpha}$ is the typical distance
that a particle travels) then $\P^{\alpha,n}[K_n]\to 1$  as $n\to \infty$
uniformly in $S \in A_n$. We have
\[
\P^{\alpha,n}[S_t \in A_n]  \le  \P^{\alpha,n}[H_n, \, K_n \, \mathrm{and}\, S_t \in A_n]  + \P^{\alpha,n}[H_n^c] + \P^{\alpha,n}[K_n^c].
\]
On the event $H_n \cap K_n$, in order for $S_t\in A_n$, all points which have been updated and have travelled
to a different $r_n$-square have to distribute themselves among the $O(n^2 /r_n^2)$ $r_n$-squares in
$[0, an] \times [0,bn]$ in such a way that none of these squares receives more than
$(\frac{n}{r_n})^\eps |S| \frac{r_n^2}{n^{2\alpha}}$ such points.
(Note here that we use the fact that in order for $S_t \in A_n$ then $S_t\subseteq [0, an] \times [0,bn]$
and thus $r_n$-squares outside of $[0, an] \times [0, bn]$ do not contribute.)
It follows that at least $(\frac{r_n}{n})^\eps |S|^{-1} \frac{n^{2\alpha}}{r_n^2}
n^{3/4-2\eps}$ $r_n$-squares have to intersect $S_t$. It remains to
choose $r\leq n^{1-\eps}$ and $r_n \ll n^\alpha$ so that this number
exceeds $(n/r_n)^{3/4+\eps}$ (recall that we already required $r_n$ to satisfy
$n^{3/4-\eps}\frac{r_n^2}{n^{2\alpha}} \ge 1$).
Let $r_n:= n^\gamma$ for some exponent $\gamma$.
We can now optimize on the value of $\gamma$ under the following constraints:

\begin{equation}
\left \lbrace \begin{array}{ll}  n^\gamma \le n^{1-\eps} \\
n^\gamma \ll n^\alpha \\
1 \le  n^{3/4-\eps} n^{2(\gamma-\alpha)}  \\
\left(\frac{r_n}{n}\right)^\eps |S|^{-1} \frac{n^{2\alpha}}{r_n^2}
n^{3/4-2\eps}  \gg n^{(1-\gamma)(3/4 + \eps)}. \\
\end{array}\right.
\end{equation}
Note that since $S \in A_n$, $|S|\le n^{3/4+\eps}$ and hence
\[
\left(\frac{n}{r_n}\right)^\eps |S| \frac{r_n^2}{n^{2\alpha}} \le  n^{\eps (1-\gamma)} n^{3/4+\eps} n^{2(\gamma-\alpha)}\,.
\]
Since $\alpha>3/8$, we write $\alpha= 3/8 + \delta/2$ with $\delta>0$.
This leads to the system:
\begin{equation}
\left \lbrace \begin{array}{ll} \gamma \le 1-\eps \\ \gamma < \alpha \\ 0 \le   2\gamma -\delta -\eps    \\
 (\gamma -4) \eps + \delta - 2\gamma  >  (1-\gamma)\eps -  3/4\gamma\end{array}\right.
\end{equation}
Since $\eps$ can be chosen arbitrarily small, we just need to find
$\gamma\in(0,1)$ so that $ 2\gamma > \delta$, $\gamma < \alpha$
and $5 \gamma /4 < \delta$. One sees that $\gamma := \frac 3 5 \delta$ does this.
It is easy to see that if $\alpha\le 3/8$ (so that $\delta \le 0$),
the above system does not have a solution.
\fbox{}

\noindent{\bf Remark:}
Note that the proof implies quantitative noise sensitivity results in that
\[
\lim_{n \to \infty}\E[f_n(\eta^{\alpha,n}_0)f_n(\eta^{\alpha,n}_{n^{-\beta}})]-\E[f_n(\eta^{\alpha,n}_0)]^2=0
\]
for some exponent $\beta>0$.
The asymptotic independence should be true for any $\alpha>0$ and $\beta<3/4$, but the present proof
would only give an exponent $\beta=\beta(\alpha)$ which would converge to zero as $\alpha \to 3/8$.

\medskip\noindent
The proof of Proposition \ref{pr.SR2} relied on results from
\cite{GPS} concerning the mesoscopic structure of $\Spec_n$.
To prove Theorem \ref{th.SR} the idea is to bootstrap
the proof of Proposition \ref{pr.SR2} using a two-scale argument.

The two scales are described as follows.
Let $\alpha$ be any positive exponent such that $\alpha<1/2$ (we already know how to deal with $\alpha>3/8$).
Let us partition our domain $[0,an]\times[0,bn]$ into a coarse grid of $R$-squares with $R:= n^{2\alpha}$
and a finer grid of $r$-squares with $r:= R^\gamma$ (where $\gamma \in (0,1)$ will be chosen later).

\begin{figure}
\begin{center}
\includegraphics[width=0.95\textwidth]{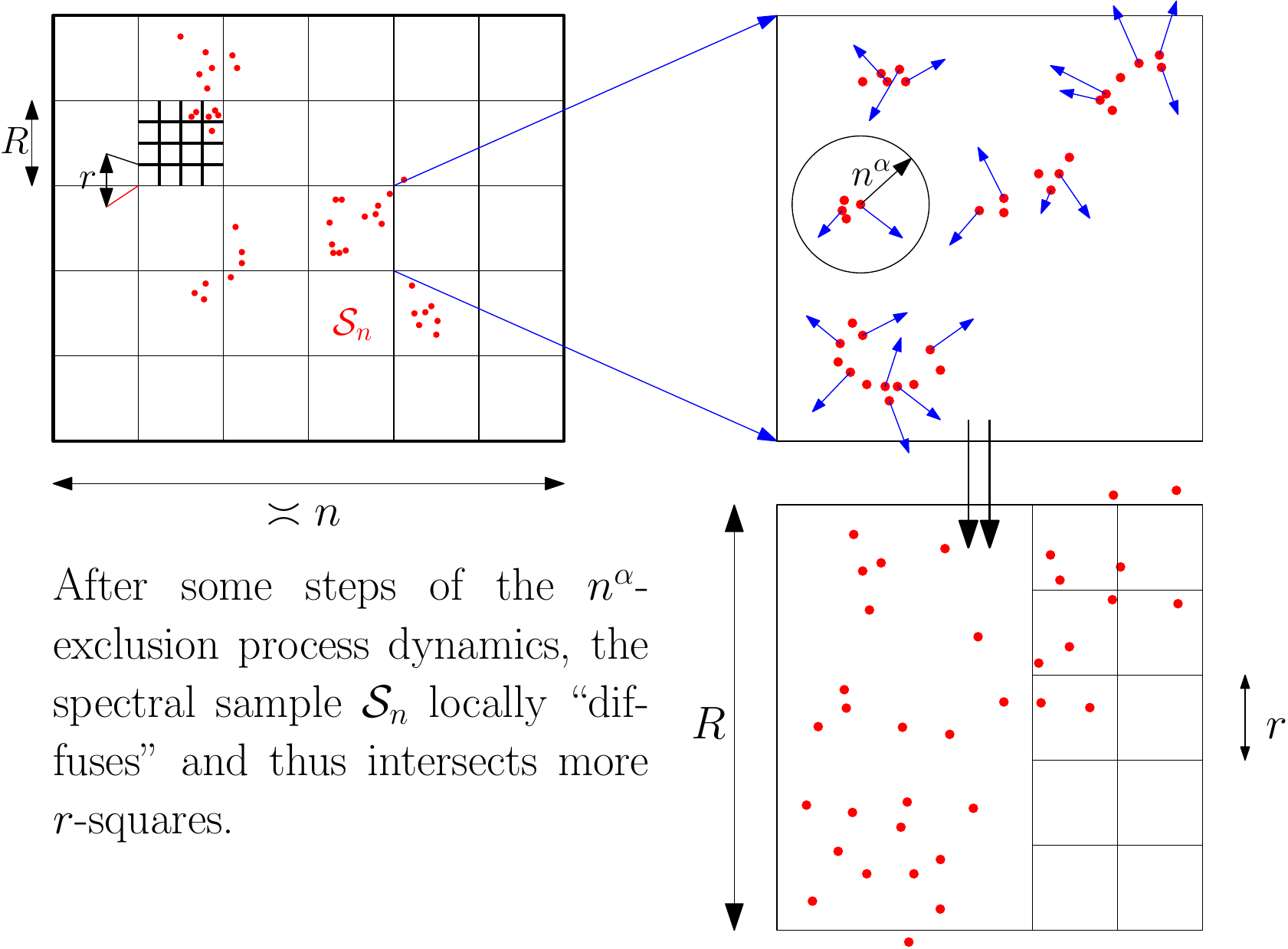}
\end{center}
\caption{\bf A visualization of the proof of Theorem \ref{th.SR}}\label{f.Bootstrap}
\end{figure}

The idea in the proof of Theorem \ref{th.SR} is that
if the spectral sample $\Spec_n$ intersects some $R$-square $Q$, then it should behave more or less like
a sample of ``$\Spec_R$'' within this square. In fact, to apply the technology of the proof of Proposition \ref{pr.SR2}
it is enough for us to prove that with probability going to one, one can find some $R$-square  $Q$,
for which $\Spec_n \cap Q$ satisfies the analogue of the event $A_n$ of Proposition \ref{pr.SR2}
but defined in terms of $Q$. Now, in $Q$ the range of the exclusion dynamics is
$n^\alpha = R^{1/2}$, so locally, relative to the size of our $R$-square, the range-exponent is larger than $3/8$.
See Figure \ref{f.Bootstrap} for a schematic view of the strategy of the proof.

To justify this intuition of local behavior we will need the following two lemmas.
Here, given an $R$-square $Q$, $2Q$ refers to the square concentric to $Q$ but with sidelength $2R$.

We now need to give a definition.
\begin{definition} \label{defn:jointpivotal}
Let $f: \{0,1\}^{[n]} \to \{-1,1\}$ and
$x_1,\ldots,x_k$ be some subset of the $n$ bits.
The event that ``$x_1,\ldots,x_k$ are \emph{jointly pivotal}'' is the event
that when the other bits are chosen, the induced Boolean function
of $x_1,\ldots,x_k$ is a function which depends on {\bf all} $k$ bits
(in the sense that for each bit $x_i,$ there is a choice of the other bits
$x_1,\ldots,x_k$ making $x_i$ pivotal).
\end{definition}

Note that the event that $x_1,\ldots,x_k$ are jointly pivotal is an event
which is measurable with respect to the complementary bits.

The following lemma follows quite easily from Lemma 2.1 in \cite{GPS}.

\begin{lemma} \label{lemma:spectrumANDjointpiv}
Let $k\geq 1$. For any  $f: \{0,1\}^{[n]} \to \{-1,1\}$ and any $x_1, \ldots , x_k$,
\[
\Pb{x_1, \ldots, x_k \in \Spec_n} \le \Pb{x_1, \ldots, x_k \text{ are jointly pivotal}}.
\]
\end{lemma}

The proof of the following lemma relies on the technology which was developed 
in \cite{PivotalMoments}, where the $k$-th moments of the set of pivotals,
$\Eb{|\mathcal{P}_n|^k}$, for percolation crossing events were studied.
This technology, rather than the actual results explicitly stated in 
\cite{PivotalMoments} are needed.  Indeed, we will explain below that the
following lemma follows from a slight strengthening of 
Theorem 1 in \cite{PivotalMoments}. It turns out that the proof of Theorem 
1 in \cite{PivotalMoments} in fact provides this strengthened statement 
we need. Since such a statement was unfortunately not explicitly 
stated in \cite{PivotalMoments}, we give some details below.  

\begin{lemma}\label{l.MomentBound}
Let $\eps>0$. For any $k\geq 1$, there exists a universal constant
$C_k=C_k(\eps)>0$ such that  for any $1\le r \le R\le n$ and any $R$-square $Q$
in $[0,an]\times [0,bn]$, if $\Spec_{(r)}^{2Q}$
denotes the set of $r$-squares in $2Q$ intersected by $\Spec_n$, then
\[
\Eb{|\Spec_{(r)}^{2Q}|^k \md \Spec_n \cap 2Q \neq \emptyset} \le C_k \Bigl[ (R/r)^{3/4+\eps/2} \Bigr]^k\,.
\]
\end{lemma}

\noindent{\bf Proof (sketch).} We assume that $a=b=1$. 
Let us start with the case $r=1$ and $R=n$, 
in which case $\Spec_{(r)} = \Spec=\Spec_n$. 
Take $k\geq 1$. Lemma \ref{lemma:spectrumANDjointpiv} yields
\begin{eqnarray*}
\lefteqn{\Eb{|\Spec_n|^k}=\sum_{x_1, \ldots, x_k}\Pb{x_1, \ldots, x_k \in \Spec_n}}\\
& & \leq \sum_{x_1, \ldots, x_k}\Pb{x_1, \ldots, x_k \text{ are jointly pivotal}}.\
\end{eqnarray*}

Now, Theorem 1 in \cite{PivotalMoments} states that if $\mathcal{P}_n$ 
denotes the set of pivotal points of our left
to right crossing event $f_n$, then for every $k\geq 1$, as $n \to \infty$, one has 
\[
\Eb{|\mathcal{P}_n|^k} = n^{3 k/4 + o(1)}\,.
\]

Since 
\[
\Eb{|\mathcal{P}_n|^k} = \sum_{x_1, \ldots, x_k} \Pb{x_1,\ldots, x_k \in \mathcal{P}_n},
\]
the above theorem from \cite{PivotalMoments} would imply our lemma if one 
had an upper bound of the form 
\[
\Pb{x_1, \ldots, x_k \text{ are jointly pivotal}} \le \Pb{x_1, \ldots, x_k \text{ are pivotal}}\,.
\]
Unfortunately, only the reversed inequality holds and 
therefore one has to look more closely into the proof in \cite{PivotalMoments}.

The strategy in \cite{PivotalMoments} in order to control $\Eb{|\mathcal{P}_n|^k}$ is to show that the ``main'' contribution in the sum 
$\sum_{x_1,\ldots, x_k} \Pb{ x_1, \ldots, x_k \text{ are pivotal} }$ comes from $k$-tuples of points $(x_1,\ldots, x_k)$ which are 
``macroscopically far apart''.
More precisely, their proof shows that for every $k$,
there exists a constant $c>0$ small enough such that 

\begin{align*}
\Eb{|\mathcal{P}_n|^k} & = \sum_{x_1,\ldots, x_k} \Pb{ x_1, \ldots, x_k \text{ are pivotal} } \\
& \le  2 \sum_{ x_1, \ldots, x_k :\, \inf_{i\neq j} |x_i - x_j| > cn} \Pb{ x_1, \ldots, x_k \text{ are pivotal} } \,.
\end{align*}

In order to achieve this, they introduce a way to organize (or classify) the possible $k$-tuples $(x_1, \ldots, x_k)$. They call this 
classification the {\it disjoint box method} which is elaborated in 
Section 4 in \cite{PivotalMoments}. In order to apply this {\it disjoint box method},
one needs the separation of arms and quasi-multiplicativity from the works of Kesten (see \cite{WWperc} for a more modern account). 
The same tools apply, when we study {\it jointly} pivotal points instead of actual pivotal points. Therefore following the {\it disjoint box method} and the subsequent sections in \cite{PivotalMoments},
one has the existence of some $\bar c>0$ so that 

\begin{align*}
 & \sum_{x_1,\ldots, x_k} \Pb{ x_1, \ldots, x_k \text{ are jointly pivotal} }  \\
 &  \hskip 3 cm \le  2 \sum_{ x_1, \ldots, x_k :\, \inf_{i\neq j} |x_i - x_j| > \bar cn} \Pb{ x_1, \ldots, x_k \text{ are jointly pivotal} } \\
 &   \hskip 3 cm  \le O(1) n^{2k} \alpha_4(n)^k \,,
\end{align*}
by independance on disjoint sets and quasi-multiplicativity. Using $\alpha_4(n) = n^{-5/4 +o(1)}$ (from \cite{SW}), this implies our result for $r=1$
and $R=n$.

Now, if $1\le r \le n$ is some intermediate scale between 1 and $n$ and we 
still have $R=n$, one needs the following extension of 
Lemma \ref{lemma:spectrumANDjointpiv}, namely that 
if $B_1, \ldots, B_k$ are any $k$ disjoint boxes, then 
\begin{align*}
& \Pb{\Spec_n \text{ intersects all the } B_i, i\in \{1,\ldots, k\} }   \\
& \hskip 3 cm \le \Pb{B_1,\ldots, B_k \text{ are jointly pivotal boxes} }\,,
\end{align*}
where the definition of jointly pivotal disjoint boxes is the straightforward generalization of Definition \ref{defn:jointpivotal}. Note that this 
inequality also follows immediately from Lemma 2.1 in \cite{GPS}. 

Let us divide the rectangle into $O(n^2/r^2)$ $r$-squares $B_i$ so that 
\[
\Eb{|\Spec_{(r)}|^k} = \sum_{B_{i_1}, \ldots, B_{i_k}} \Pb{\Spec_n \text{ intersects all the } B_{i_l}, l\in \{1,\ldots, k\}}
\]

Thanks to quasi-multiplicativity which in a way factorizes what happens above and below scale $r$, the same {\it disjoint box method} from Section 4
in \cite{PivotalMoments} applies in the present setting and yields
\begin{align*}
& \sum_{B_{i_1}, \ldots, B_{i_k}} \Pb{B_{i_1},\ldots, B_{i_k} \text{ are jointly pivotal boxes }} \\
& \hskip 3 cm \le  O(1) \left(\frac {n^2}{r^2}\right)^k \alpha_4(r,n)^k\,,
\end{align*}
which implies the desired result.
For general values of $1\leq r\leq R$, similar considerations apply to the moments of 
$|\Spec_{(r)}^{2Q}|$ conditioned on $\Spec_n$ intersecting $2Q$.
\fbox{}
\vskip 0.3 cm

The second lemma, which says that points do not travel far, is the following.
\begin{lemma}\label{l.travel}
Let $\tilde S$ be any subset of the $R$-square $Q$, and let $\tilde S_t=\pi^{n,n^\alpha}_t( \tilde S)$ where
$n^\alpha$ refers to the $n^\alpha$-exclusion process dynamics. Then for some $c=c(t)>0$,
\[
\P(\tilde S_t \subseteq 2Q)\geq 1-R^2e^{-c R}.
\]
\end{lemma}
\noindent{\bf Proof of Lemma \ref{l.travel}.}
Each point in $x\in \tilde S$ performs a simple random walk $(X_t)_{t\geq 0}$,
where at rate one $X_t$ jumps to a uniform point in $B(X_t, n^\alpha)$.
Therefore the variance at time $t$ of the random walk is
$\asymp t n^{2\alpha} =t\, R$. In addition, by standard arguments about random walks,
there is some $c=c(t)>0$ such that $\Pb{|X_0 - X_t| > R} \le e^{-c R}$.
Hence,
\[
\P(\tilde S_t \not \subseteq 2Q)\leq |\tilde S|e^{-c R}\leq R^2e^{-c R}.
\]
\fbox{}\\

\medskip\noindent
We now have all the ingredients that we need, we refer the reader to Figure \ref{f.Bootstrap}
for the idea of the proof.\\
\noindent{\bf Proof of Theorem \ref{th.SR}.}

Recall that $R^{1/2}=n^\alpha,$ $r=R^\gamma$ and let $\eps,\gamma>0$ be
chosen later. Define $A_n= A^1_n \cap A^2_n$ where


\begin{eqnarray*}
& & A^1_n=\left\{
\begin{array}{l}
S\subseteq [0,an]\times[0,bn]\cap \T: \text{there exists an $R$-square $Q$ in our} \\
\textrm{ coarse grid s.t. } R^{3/4 -\eps} \le |S\cap Q| \le R^{3/4+\eps}
\end{array}
\right\} 
\\
& & A^2_n=\left\{
\begin{array}{l}
S\subseteq [0,an]\times[0,bn]\cap \T:\text{for all $R$-squares $Q$ }  \\
\textrm{ in the coarse grid, the number of $r$-squares inside } \\
\textrm{ $2Q$ which intersect $S$ is less
than $\left(\frac R r \right)^{3/4+\eps}$}
\end{array}
\right\}
\end{eqnarray*}


It follows directly from \cite{GPS} that, conditioned on $\Spec_n \neq \emptyset$,
with conditional probability tending to one, there exists
at least one $R$-square $Q$, such that $|\Spec_n \cap Q| \geq R^{3/4-\eps}$.
(This is in fact the key thing which makes the main proof work in
\cite{GPS}). Now to prove that $\Pb{A^1_n \md \Spec_n \neq \emptyset} \to 1$,
we will show that the square we found does not
contain too many spectral points. For this we shall prove the sharper result:
\begin{equation} \label{eqn8}
\Pb{\text{For all $R$-squares $Q$ in the coarse grid, } |\Spec_n \cap Q| \le R^{3/4+\eps} }  \to 1.
\end{equation}
Assuming this, we obtain that $\Pb{A^1_n \md \Spec_n \neq \emptyset} \to 1$.
Showing (\ref{eqn8}) basically corresponds to showing $\Pb{A^2_n} \to 1$ had $r$ been 1.
We prove that $\Pb{A^2_n} \to 1$; (\ref{eqn8}) is proved in the same way.

From Lemma \ref{l.MomentBound}, we find that
\begin{eqnarray*}
\lefteqn{\Pb{ \text{There is an $R$-square $Q$, s.t.
$|\Spec_{(r)}^{2Q}| > (R/r)^{3/4+\eps}$} }}\\
& & \le \;  \sum_{\text{squares $Q$}} \Pb{|\Spec_{(r)}^{2Q}| > (R/r)^{3/4+\eps}} \\
& & =  \sum_{\text{squares $Q$}} \Pb{\Spec_n \cap 2Q \neq \emptyset }
 \Pb{|\Spec_{(r)}^{2Q}| > (R/r)^{3/4+\eps} \md \Spec_n \cap 2Q \neq \emptyset} \\
& &  \le \;   C_k \, \Bigl( \frac r R \Bigr)^{k \eps/ 2} \sum_{\text{squares $Q$}} \Pb{\Spec_n \cap 2Q \neq \emptyset } \\
& &  \le \; C_k \, \Bigl( \frac r R \Bigr)^{k \eps/ 2}9\Eb{|\Spec_{(R)}|} \\
& &  \le C_k \, \Bigl( \frac r R \Bigr)^{k \eps/ 2} 9\Bigl( \frac n R\Bigr)^{3/4 +\eps} \\
& &  = 9C_k\, \Bigl( R^{\gamma - 1}\Bigr)^{k \eps/ 2}
      \Bigl( \frac n {n^{2\alpha}}\Bigr)^{3/4 +\eps} \; \text{since we fixed }\left\lbrace
    \begin{array}{l} R^{1/2}= n^{\alpha}\\ r = R^\gamma. \end{array} \right.
\end{eqnarray*}
The term 9 arises since we can cover any fixed $2Q$ box by 9 $Q$-boxes and the last inequality
follows (as before) by \cite{GPS}.

Therefore, for any arbitrary choice of $\gamma, \eps \in (0,1)$, we can
choose $k=k(\gamma, \eps)$, so that this probability goes to zero
as $n\to\infty$. Hence, $\P(A^2_n \md \Spec_n \neq \emptyset) \to 1,$ and so we conclude that
$\P(A_n \md \Spec_n \neq \emptyset) = \P(A_n^1 \cap A_n^2 \md \Spec_n \neq \emptyset) \to 1$ as $n\to \infty$.
\vskip 0.3 cm

The rest works exactly as in the proof of Proposition \ref{pr.SR2}: for any set $S\in A_n$ with
$A_n=A_n^1 \cap A_n^2$, one chooses an $R$-square $Q=Q(S)$ so that
$R^{3/4-\eps} \le |S\cap Q| \le R^{3/4+\eps}$ (there exists such a square since $S\in A_n^1$).
We thus have a ``large'' set  $\tilde S := S^Q$ which by Lemma \ref{l.travel} is very likely to stay in $2Q$,
and by the same argument as in Proposition
\ref{pr.SR2}, one can show that if one chooses $\gamma\in (1/8,1/5)$
and $\eps$ small enough, then with
high probability $\tilde S_t$ will intersect many $r$-squares
in $2Q$, which because of the event $A_n^2$ will force $S_t \supseteq \tilde S_t$ to leave the set $A_n$.
We then apply Theorem \ref{thm3}.\fbox{}\\

\medskip\noindent
{\bf Remarks:} 
1. One could hope (in the spirit of  the remark after the proof of Proposition \ref{pr.SR2}) that
the above bootstrap procedure would not only give a better control on
the exponent $\alpha$ but would also imply better bounds on the exponent $\beta$ (see the remark
after the proof of Proposition \ref{pr.SR2}
for the notation). This is
true in a way: one can check that the above proof provides a positive exponent $\beta$ for all $\alpha >0$.
Nevertheless with this type of technique, $\beta$
would also have to go to zero as $\alpha\to 0$. \\
2. It might also seem from this ``multi-scale'' proof, that one could
push the argument all the way to very local dynamics  (say, in $\log n$
instead of $n^\alpha$ and maybe even all the way to the elegant
nearest-neighbor dynamics). However, if one looks into the proof of
Theorem \ref{th.SR}, there is at least one place which is quite
problematic: the proof relies on the moment bounds given by Lemma
\ref{l.MomentBound}. If, for example, our dynamics was a
$\log n$-medium range dynamics, then $\frac r R \geq (\log n)^{-1}$,
and one would need very sharp upper tail estimates on
$|\Spec_{(r)}^{2Q}|$ to be able
to conclude the proof. In fact by looking into what should be the fluctuation
of the random variable $|\Spec_{(r)}^{2Q}|$,
one can see that there is no hope that
the technique we are presently using could handle the case of the
nearest-neighbor dynamics.

\section{Discussion on sharper results concerning Exclusion Sensitivity for percolation}
\label{secdiscussions}

\subsection{Quantitative Exclusion sensitivity}
Corollary \ref{corr1} shows that if one has quantitative bounds on the
noise sensitivity of a sequence of monotone Boolean functions, then the
same quantitative bounds hold for Complete Graph exclusion sensitivity.
It is not hard to check that using our spectral approach (Section
\ref{secthm5}), one can also obtain quantitative bounds on
exclusion sensitivity even when the underlying graph is not the
complete graph. See the remark after the proof of Proposition \ref{pr.SR2}.

Nevertheless, for the medium-range $n^\alpha$ dynamics, as $\alpha$ goes to zero,
the polynomial bounds on the exclusion sensitivity that follow from our ``spectral
proof'' become weaker and weaker.
In contrast one would believe that the true bound should remain of the same order as
in the i.i.d.\ dynamics (i.e. the threshold of exclusion sensitivity
should happen near $\eps_n \asymp t_n \approx n^{-3/4}$ independently of the
medium-range exponent $\alpha$).
Such a sharp control on the exclusion sensitivity seems to be quite challenging with
the present techniques in this paper. 

\subsection{Medium-range dynamics: coupling with i.i.d. dynamics}
In Section \ref{s.shortrange}, we used our spectral approach to obtain exclusion sensitivity of percolation in the case of
medium-range dynamics. Instead we could 
have used the approach of Section \ref{sect:monotonicity}.
Let us briefly sketch how this would work. If one considers an $n^\alpha$-medium-range dynamics in $[0,n]^2$, then divide $[0,n]^2$
into $n^{2(1-\alpha)}$ squares of sidelength $n^\alpha$. In each of these squares, one can couple the medium-range dynamics with an i.i.d. dynamics
up to a fluctuation of $O(n^\alpha)$ bits. (Furthermore, there are some boundary issues to deal with since
neighboring boxes interact with each other, but let us neglect this here.)
In total, this gives a fluctuation of $O(n^{2(1-\alpha)} n^{\alpha}) = O(n^{2-\alpha})$ bits between an i.i.d.
dynamics and a $n^\alpha$-medium-range dynamics.
From the known behavior of \emph{near-critical percolation} on the triangular lattice, it follows that if one
adds $n^{2-\alpha}$ random points to a uniform configuration, and if
$n^{2-\alpha}\ll n^{5/4}$, then it is very unlikely that there is a change in the crossing event $C_n$.
Using this information, one can potentially obtain
medium-range exclusion sensitivity with this coupling method
\emph{but only if the exponent $\alpha$ satisfies $\alpha>3/4$}. This shows that ``paradoxically'',
the coupling method seems weaker than the
spectral approach for the study of medium-range dynamics but on the other hand, it seems more robust for proving quantitative exclusion sensitivity results.
Note in particular that there is no hope to study the nearest-neighbor dynamics simply via a coupling with the i.i.d.\ dynamics.

\subsection{Limitations of our spectral approach for medium-range dynamics and nearest-neighbor dynamics}

In Section \ref{s.shortrange}, in order to detect a singular behavior of $S_t$,
when $S\sim \nu_{C_n}$ (where $\nu_{C_n}$ denotes the spectral measure of $I_{C_n}$),
we relied essentially on the sizes of $\Spec_n$ or $\Spec_{(r)}$ restricted to an $R$-mesoscopic box $Q$, but not so much
on the local geometry of the spectrum within $Q$. Even though this would be a tempting way to obtain better estimates, there is a reason why
we stuck to ``numerical'' quantities: this follows from the fact that the technology introduced in \cite{GPS}
is at the moment not robust enough to imply \emph{0-1 laws}
for local geometric schemes.



Let us now illustrate why it would be very difficult to apply the present spectral approach in the case of nearest-neighbor dynamics.
We believe the following result should hold.
\begin{conjecture}
If $\{\eta_t\}$ denotes the {\bf nearest-neighbor} dynamics on the triangular lattice (at $p_c=1/2$), then
\[
\lim_{n\to \infty}\E[f_n(\eta_0) f_n(\eta_{t_n})]- \E[f_n(\eta_0)]^2=0,
\]
as long as $t_n \gg (n^2 \alpha_4(n))^{-1} \approx n^{-3/4}$.
\end{conjecture}

If one wanted to use the spectral
approach developed in this paper to prove this result, one would need to
find a property (an event $A_n$) almost surely satisfied for $\Spec_n$,
but almost never satisfied once a tiny portion of the $\approx n^{3/4}$
points in $\Spec_n$ have moved roughly one or two (neighboring) steps!
Consider the following analogue: start a random walk on $\Z$ at
$X_0:= \sqrt{n}$ and consider the set
$S:= \{ k\in[0,n] \text{ s.t. } X_k=0  \} $ ($S$ could be empty).
Now, perturb $S$ slightly into $\tilde S$ by choosing $\log n$ points at
random in $S$ (if $S$ is non-empty then it is of size $O(\sqrt{n})$) and
let these points move by 2 or 4 steps (with the exclusion rule), where 2
or 4 is chosen in order to preserve the parity constraints. Can one tell
whether $\tilde S$ looks significantly different from $S$?
I.e. are $\tilde S$ and $S$ singular?

\medskip \noindent
We note that it could be that when the dynamics are too local, then the spectral approach
fails, since with this approach we do not take into account the
possible sign-cancellations occurring in \eqref{eqn1}.

\subsection{Exceptional times}

The complete-graph dynamics is certainly not suitable for a full-plane study,
but our medium-range dynamics is. One might ask whether there are exceptional
times for these conservative dynamics at which infinite clusters appear.
Such questions are asked in \cite{SS, GPS} for the SNPS process.
The trouble with our present techniques is that if one fixes a
medium-range dynamics with $L:=n^{\alpha}$,
then for scales much higher than $L$, the dynamics would start looking so
local that our bounds would not be good enough to imply existence of
exceptional times.


\section*{Acknowledgements}

We would like to thank Cedric Bernardin and Milton Jara for enlightening
discussions about some FKG-type properties for Exclusion processes.

We would also like to thank the Mittag-Leffler Institute for their support.
It was during a visit there that this project was initiated and a substantial
part of the work was performed.

The research of C.G.\ was partially supported by ANR grant
BLAN06-3-134462. The research of J.E.S.\ was
supported by the Swedish Natural Science Research Council
and the G\"{o}ran Gustafsson Foundation for Research in
Natural Sciences and Medicine. The research of E.I.B.\ was also
supported by the G\"{o}ran Gustafsson Foundation for Research in
Natural Sciences and Medicine.

Finally, we would like to thank the anonymous referee who read the
manuscript very carefully and provided us with a large number of very 
useful comments.

\noindent {\bf Erik Broman}\\
Mathematical Sciences \\
Chalmers University of Technology \\
and \\
Mathematical Sciences \\
G\"{o}teborg University \\
SE-41296 Gothenburg, Sweden \\
{broman@chalmers.se} \\
{\tt {http://www.math.chalmers.se/\string~broman/}}

\medskip\noindent {\bf Christophe Garban}\\
CNRS \\
Ecole Normale Sup\'erieure de Lyon, UMPA \\
46 all\'ee d'Italie \\
69364 Lyon Cedex 07 France \\
{\tt {http://www.umpa.ens-lyon.fr/\string~cgarban/}}.

\medskip\noindent {\bf Jeffrey E. Steif}\\
Mathematical Sciences \\
Chalmers University of Technology \\
and \\
Mathematical Sciences \\
G\"{o}teborg University \\
SE-41296 Gothenburg, Sweden \\
{steif@math.chalmers.se} \\
{\tt {http://www.math.chalmers.se/\string~steif/}}

\end{document}